\pdfoutput=1
\documentclass[3p]{elsarticle}
\makeatletter
\def\ps@pprintTitle{%
 \let\@oddhead\@empty
 \let\@evenhead\@empty
 \def\@oddfoot{\centerline{\thepage}}%
 \let\@evenfoot\@oddfoot}
\makeatother
\usepackage{amsmath}
\usepackage{amssymb}
\usepackage{hyperref}
\hypersetup{pdfauthor={Name}}
\usepackage{xcolor}
\usepackage{epstopdf}
\usepackage[T1]{fontenc}
\usepackage{subcaption}

\newcommand{\tf}{{\rm f}}
\newcommand{\R}{\mathbb{R}}

\begin{document}

\begin{frontmatter}

\title{Intrinsic Integration\tnoteref{funding}}
\tnotetext[funding]{This work was funded by the National Institutes of Health under grant no.~R01-HL-143350 and Army Research Office under grant no.~W911NF-18-1-0281.}

\author[1]{Navdeep Dahiya}
\ead{ndahiya3@gatech.edu}
\author[2]{Martin Mueller}
\ead{martin.mueller@gatech.edu}
\author[1]{Anthony Yezzi\corref{cor1}}
\ead{ayezzi@ece.gatech.edu}

\cortext[cor1]{Corresponding author}
\address[1]{Electrical \& Computer Engineering, Georgia Institute of Technology, Atlanta, GA, USA}
\address[2]{Rivian Automotive, Palo Alto, CA, USA}


%
%

\begin{abstract}
 If we wish to integrate a function $h|\Omega\subset\Re^{n}\to\Re$ along
a single $T$-level surface of a function $\psi |\Omega\subset\Re^{n}\to\Re$, then a number
of different methods for extracting finite elements appropriate to
the dimension of the level surface may be employed to obtain an explicit representation over which
the integration may be performed using standard numerical quadrature
techniques along each element. However, when the goal is to compute
an entire continuous family $m(T)$ of integrals over all the
$T$-level surfaces of $\psi$, then this method of explicit level set extraction is no longer practical.
We introduce a novel method to perform this type of numerical integration efficiently by making
use of the coarea formula. We present the technique for discretization of the coarea formula and present the
algorithms to compute the integrals over families of T-level surfaces. While validation of our method in the special case of a single level surface demonstrates accuracies close to more explicit isosurface integration methods, we show a sizable boost in computational efficiency in the case of multiple T-level surfaces, where our coupled integration algorithms significantly outperform sequential one-at-a-time application of explicit methods.
\end{abstract}


\end{frontmatter}


\section{Introduction}
Level set methods for capturing moving fronts introduced by Osher and Sethian~\cite{OsherSethian87} have proven to be a robust numerical device for large collection of diverse problems. These methods are frequently used for propagating interfaces in computational physics~\cite{OsherFedkiw03,Sethian99} and fluid dynamics~\cite{SethianSmerka03}, constrained optimization~\cite{OsherSantosa01,OsherSantosa08}, and computer vision~\cite{Cremers06,tsaiyezzi} among others. Consequently, the problem of numerical integration of a function over volumes and surfaces implicitly defined by a fixed iso-contour/surface of a level set function frequently arises. Over the years many different approaches and methods have been developed for these types of numerical integration problems.

The traditional approach for evaluating such integrals has been to make use of the regularized Dirac delta or the Heaviside function~\cite{OsherSethian87}. In~\cite{TornbergEnquist03}, the authors showed the standard regularization approaches to be inconsistent and that they may lead to non convergent solutions with O(1) errors. Enquist et al.~\cite{TornbergEnquist05} later presented discretization techniques for the Dirac delta that were first order accurate with a second order convergence rate. In~\cite{Smereka06}, Smereka further developed a second order accurate discretization of the Dirac delta function with application to calculating line and surface integrals in level set based methods.

Another class of methods explicitly reconstruct the interface to find a faceted mesh representation of the surface through various methods like marching cubes~\cite{LorensenMarchingCubes} or marching tetrahedra~\cite{LingChen98,Hummel95} and then apply standard quadrature schemes to the surface and volume elements. In~\cite{MinGibou07}, Min and Gibou presented a geometric technique for calculating such integrals involving discretizing the interface into simplices and using numerical integration quadrature rules on the simplices. They showed the approach to be second-order accurate and being robust to interface perturbation on the grid. They showed several practical examples and compared their results to the techniques introduced by Smereka in~\cite{Smereka06}. 

In yet another approach making use of the well known divergence theorem, M\"uller et al.~\cite{Mueller12} presented a simple tree based approach and showed it to be comparable to the methods of Min-Gibou and Smereka in accuracy and convergence. In~\cite{Mueller13} M\"uller et al. presented a new method based on the solution of a small linear system based on a simplified variant of moment-fitting equations. The authors showed the method to be several orders of magnitude more accurate than earlier approaches. More recently,~\cite{Saye15} presented a high-order accurate numerical quadrature method for evaluating integrals over curved surfaces and volumes defined via a fixed isosurface of a function restricted to a given hyperrectangle. The method converts the implicitly defined geometry into the graph of a height function leading to a recursive algorithm on the number of dimensions requiring one-dimensional root finding and one-dimensional Gaussian quadrature. 

In~\cite{kublik2013implicit}, the authors proposed a formulation for computing integrals of the form, $\int_{\partial \Omega} v\left(\mathbf{x} \left(s\right)\right)ds$ in the levelset framework, when the domain $\Omega$ is represented implicitly by the signed distance function to its boundary $\partial \Omega$. The authors expressed this integral as an average of integrals over nearby level sets of $d_{_{\partial \Omega}}$ , where
these nearby level sets continuously sweep a thin tubular neighborhood around the boundary ${\partial \Omega}$ of radius $\epsilon$. With this formulation, the authors proposed a numerical method based on integral equation formulations for solving the Poisson problem with constant coefficients, subject to Dirichlet, Neumann, Robin or mixed boundary conditions. This proposed formulation involved an exact formulation for computing boundary integrals in the level set framework, and provided a natural way of defining and computing boundary integrals in applications using the closest point formulations~\cite{RuuthClosestPoint1,RuuthClosestPoint2}. In~\cite{kublik2016integration}, the authors further extended this technique to include integration along curves in three dimensions. Here the authors proposed a new formulation using the closest point mapping for integrating over smooth curves and surfaces with boundaries that are described by their closest point mappings. The authors studied various aspects of this formulation and provided a geometric interpretation in terms of the singular values of the Jacobian matrix of the closest point mapping. In~\cite{kublik2018extrapolative}, the authors further extended the framework of~\cite{kublik2013implicit, kublik2016integration}, and proposed an extrapolative approach for computing integrals over a class of piecewise smooth hypersurfaces, given implicitly via a level set function. Although this method was based on the classical
approximation used in the level set framework that smears out the Dirac $\delta$-function to a bump function with a compact support, in this work the authors combined the classical formulas for the more challenging cases in which the hypersurfaces have kinks and corners. The authors made of use
special kernels with vanishing moments for the approximation of the Dirac $\delta$-function, and their proposed method did not require any local explicit parameterization of the hypersurface, nor the explicit locations of the corners and kinks on the hypersurface.

As we have seen so far that if we wish to integrate a function $h|\Omega\subset\Re^{n}\to\Re$ along
a single $T$-level surface of a function $\psi$$ $, then a number
of different methods for extracting finite elements appropriate to
the dimension of the level surface (line segments, triangles, tetrahedra,
etc.) may be employed to obtain an explicit representation over which
the integration may be performed using standard numerical quadrature
techniques along each element. However, when the goal is to compute
the an entire continuous family $m(T)$ of integrals over all the
$T$-level surfaces of $\psi$, 
\[
m(T)=\int_{\psi^{-1}(T)}h\,ds
\]
then this method of explicit level set extraction is no longer practical.
An alternative is to first construct the anti-derivative $M(T)$
\[
M(T)=\int_{-\infty}^{T}\int_{\psi^{-1}(\tau)}h\,ds\,d\tau=\int_{\psi(x)<T}h\|\nabla\psi\|\,dx
\]
which, by the coarea formula, is obtained by volumetric integration
of $h\|\nabla\psi\|$ over the progressively expanding volume where
$\psi(x)\!<\!T$ (the $T$-sublevel set of $\psi$), and then differentiate
the result to obtain $m(T)=\dot{M}(T)$. 

In this paper we present the discretization of this coarea formula, derive the expressions for calculation the partial volumes across
grid cells and present an algorithm for efficiently computing $M(T)$. We also present an even more efficient algorithm for the case
where the level set function is a signed distance function. By effectively coupling the computation of the integrals over whole range of different isosurfaces, our methods lead to efficient computation of the entire family of integrals. However, it is still instructive to compare the accuracy of our methods by computing the integral for a single level surface and comparing the results to traditional single isosurface methods. Accordingly, we present results showing comparable accuracy with the methods developed in Min-Gibou~\cite{MinGibou07}. Additionally and more importantly, we can use a brute force approach of employing traditional single isosurface integration methods sequentially, one at a time, to each level surface contained in a range and compare the resulting integrals to our approach of computing the whole range of integrals in a coupled manner. We present experiments illustrating the significant computational efficiency afforded by our coupled integration methods over such a brute force approach required by traditional single isosurface methods.

The paper is organized as follows. Our main derivations and algorithms are in~\ref{sec:integrationalonglevelsets}, experimental results are in~\ref{sec:experiments},
and the conclusions follow in~\ref{sec:conclusions}.

\section{Numerical Integration along Level Sets}
\label{sec:integrationalonglevelsets}

Our goal is to compute the entire continuous family of integrals $m(T)$ of a general function $h|\Omega\subset\Re^{n}\to\Re$
over all the
$T$-level surfaces of $\psi$, 
\[
m(T)=\int_{\psi^{-1}(T)}h\,ds
\]
In order to avoid the inefficient explicit extraction of the $T$-level surfaces of $\psi$,
we first construct the anti-derivative $M(T)$
\[
M(T)=\int_{-\infty}^{T}\int_{\psi^{-1}(\tau)}h\,ds\,d\tau=\int_{\psi(x)<T}h\|\nabla\psi\|\,dx
\]
which, by the coarea formula, is obtained by volumetric integration
of $h\|\nabla\psi\|$ over the progressively expanding volume where
$\psi(x)\!<\!T$ (the $T$-sublevel set of $\psi$), and then differentiate
the result to obtain $m(T)=\dot{M}(T)$. 

\subsection{Discretization}

We will assume that the functions $h$ and $\psi$ are both sampled
over matching Cartesian grids with sample locations denoted by $x_{i}\in\Re^{n}$
and with uniform spatial increments, represented by the constant $n$-dimensional
\emph{space step vector} $\Delta x$, between samples along each orthogonal
direction\footnote{The overloaded use of $\Delta$ to represent a either a discrete interval
	size or the Lapacian operator is disambiguated by the type of symbol
	it precedes (a variable or measurement versus a function).}.
Associated to each \emph{grid point} $x_{i}$ will be the corresponding
\emph{grid cell}, denoted by ${\rm cell}(\Delta x,x_{i})$, which
refers to the $n$-dimensional rectangle centered around the grid
point $x_{i}$ and with edge lengths in orthogonal directions matching
the elements of the space step vector $\Delta x$. Finally, we will
also need to choose a \emph{range interval} $[T_{\min},T_{\max}]$
as well as a constant \emph{range increment} $\triangle T$ for our
discretized representation of the continuous functions $M$ and $m$.
Note that, unlike the space step vector $\Delta x$, the range increment
$\Delta T$ is a scalar constant.

\subsection{Cell-wise linear approximation\label{sec:linear}}

We assume that, in addition to the available sample $\psi(x_{i})$
of the level set function $\psi$ at each grid point $x_{i}$, we
have (or can compute) an estimate of the gradient $\nabla\psi(x_{i})$
at each grid point as well. We may then construct a linear approximation
$\psi_{i}$ of the level set function within each ${\rm cell}(\Delta x,x_{i})$
that matches $\psi$ and its first order derivatives at the grid point
$x_{i}$ in the center of the cell. 
\[
\psi_{i}(x)=\psi(x_{i})+(x-x_{i})\cdot\nabla\psi(x_{i}),\quad x\in{\rm cell}(\Delta x,x_{i})
\]
The minimum and maximum values of $\psi_{i}$ are attained at opposite
vertices of its corresponding cell and easily shown to be 
\begin{eqnarray*}
	\min\psi_{i} & = & \psi(x_{i})-\frac{\|\Delta x\odot\nabla\psi(x_{i})\|_{1}}{2}\\
	& = & \psi(x_{i})-\frac{\Delta x}{2}\cdot\left(\nabla\psi(x_{i})\odot\mbox{sgn}\nabla\psi(x_{i})\right)\\
	\max\psi_{i} & = & \psi(x_{i})+\frac{\|\Delta x\odot\nabla\psi(x_{i})\|_{1}}{2}\\
	& = & \psi(x_{i})+\frac{\Delta x}{2}\cdot\left(\nabla\psi(x_{i})\odot\mbox{sgn}\nabla\psi(x_{i})\right)
\end{eqnarray*}
where $\|\bullet\|_{1}$ denotes the $L^{1}$ vector norm (sum of
the absolute values of the elements), where $\odot$ denotes the Hadamard
product (element-wise multiplication), and where we have generalized
the signum function $\mbox{sgn}(\bullet)$ to act on a vector of real
numbers by simple element-wise application of the standard signum
function (which maps positive numbers to 1, negative numbers to -1,
and 0 to 0). Cell locations (along the boundary if $\nabla\psi(x_{i})\ne0)$
where these extrema are attained can be calculated using the following
formulas
\[
\arg\min\psi_{i}(x)=x_{i}-\frac{\Delta x\odot\mbox{sgn}\nabla\psi(x_{i})}{2}
\]
\[
\arg\max\psi_{i}(x)=x_{i}+\frac{\Delta x\odot\mbox{sgn}\nabla\psi(x_{i})}{2}
\]
which yield opposite vertices of the cell if and only if these extrema
are attained at unique locations. If one or more of the elements of
sgn$\nabla\psi(x_{i})$ are zero (due to the corresponding element
of $\nabla\psi(x_{i})$ being zero), then the minimum, as well as
the maximum, will be attained at multiple cell boundary locations
(including non-vertex points) rather than just at one unique vertex.
In such cases, multiple pairs of opposite vertices will exhibit the
same extremal values of $\psi_{i}$ within the cell, but the formulas
above will not yield any such pair. Such extremal vertex pair combinations
can be obtained, however, by substituting all of the zero-elements
of sgn$\nabla\psi(x_{i})$ in the formulas for $\arg\min\psi_{i}$
and $\arg\max\psi_{i}$ above with any combination of $\pm1$ (each
combination will yield a different pair of opposite min/max vertices
so long as at least one of original elements was non-zero to begin
with). If only one such vertex pair is needed, then a simple and systematic
method for obtaining one would be to replace the vector element-wise
signum function in the previous formulas with a vector element-wise
sign function
\[
S_{i}=\mbox{sign}\nabla\psi(x_{i})
\]
which produces a binary valued output vector $S_{i}$ whose elements
are either -1 (for negative input elements) or 1 otherwise (for non-negative
input elements). Using this convention, we always obtain unique set
of opposite extremal vertices ${\rm v}_{\min,i}$ and ${\rm v}_{\max,i}$
via the following formulas.
\[
{\rm v}_{\min,i}=x_{i}-\frac{\Delta x\odot S_{i}}{2}
\]
\[
{\rm v}_{\max,i}=x_{i}+\frac{\Delta x\odot S_{i}}{2}
\]
We will refer to this unique choice of opposite extremal verteces
as the \emph{minimal vertex} and the \emph{maximal vertex}, respectively,
of the cell centered around grid point $x_{i}$.

Finally, note that moving from the minimal to the maximal vertex requires
a combination of displacements along each orthogonal grid direction
with distances given by the corresponding grid step sizes stored in
the vector $\Delta x$. Since $\psi_{i}$ changes linearly, its total
increase as we move from ${\rm v}_{\min,i}$ to ${\rm v}_{\max,i}$
can likewise be decomposed as a combination of \emph{directional increases}
across the cell along each of these orthogonal grid directions, which
we represent by the vector $D_{i}$
\[
\underbrace{D_{i}}_{\mbox{increases}}=\underbrace{\Delta x}_{\mbox{distances}}\odot\underbrace{S_{i}\odot\nabla\psi(x_{i})}_{\mbox{rates of increase}}
\]
The sum of these directional increases yields the total cell increase
which can be expressed in several ways:
\[
\mbox{total cell increase}=\max\psi_{i}-\min\psi_{i}=\|D_{i}\|_{1}=\|\Delta x\odot\nabla\psi(x_{i})\|_{1}=\left(S_{i}\odot\Delta x\right)\cdot\nabla\psi(x_{i})
\]

\subsection{Cell coordinates}

Let us introduce local cell coordinates\footnote{Notice that a subscript on the cell coordinate symbol $z$ refers
	to one of its vector components rather than one of the enumerated
	grid points as for subscripts on other symbols such as $x$, $S$,
	and $\psi$. Thus $z_{k}$ denotes the $k$'th element of $z$, whereas
	$x_{i}$ denotes the location of the grid point indexed by $i$, just
	as $\psi_{i}$ denotes the linear approximation of $\psi$ over the
	corresponding grid cell indexed by $i$. Since the cell coordinates
	are local, and therefore depend upon a already-specified choice of
	grid cell, there is no need to continue representing the grid point
	index in the coordinate's notation, and so we exploit the use of subscripts
	for this other purpose here. To reinforce this difference in subscript
	meaning and reduce any resulting confusion, we reserve the use of
	the subscript symbol $i$ exclusively in reference to a grid point
	and utilize a different subscript symbol (such as $k$) whenever indicating
	a vector component. } $z=(z_{1},\ldotp,z_{n})$ which represent a shifted, scaled, and
mirrored version of the global spatial coordinates $x$, such that
a given grid cell, ${\rm cell}(\Delta x,x_{i})$ , may be locally
parameterized by the unit $n$-dimensional cube with minimal vertex
at the origin $z=(0,\ldots,0)$ and the maximal vertex at $z=(1,\ldots,1)$.
The resulting change of coordinates is expressed by the following
formulas (where we use the symbol $\oslash$ to denote vector element-wise
division or, stated more technically, multiplication by the Hadamard
inverse of the right-hand operand).
\[
x={\rm v}_{\min,i}+(\Delta x\odot S_{i})\odot z\quad\mbox{or}\quad z=(x-{\rm v}_{\min,i})\oslash(\Delta x\odot S_{i})
\]
Note that ${\rm v}_{\min,i}$ constitutes the translation, $\Delta x$
the rescaling, and $S_{i}$ the reflection. The Jacobian and its determinant
for this change of coordinates are given as follows 
\[
\frac{dx}{dz}=\mbox{diag}(\Delta x\odot S_{i}),\quad\det\left(\frac{dx}{dz}\right)=\pm\prod\Delta x
\]
where $\mbox{diag(\ensuremath{\bullet})}$ denotes the diagnonal $n\times n$
matrix whose diagonal entries match the elements of the $n$-dimensional
input vector, and where, using a slight abuse of notation, we have
written $\prod\Delta x$ to denote the product of the elements of
the vector $\Delta x$ of grid space steps. Note that the sign of
the determinant will be negative if there is an odd number of -1 entries
in the gradient sign vector $S_{i}$, otherwise the sign will be positive.
The sign is irrelevant, though, when utilizing this change of variables
in the context of integration, where only the absolute value of the
determinant is needed. 
\[
dx=\left(\prod\Delta x\right)\,dz
\]
Finally, when changing variables from $x$ to $z$ we will find the
following relation useful.
\begin{eqnarray*}
	\psi_{i}(x) & = & \psi_{i}\left(v_{\min,i}+(\Delta x\odot S_{i})\odot z\right)\\
	& = & \psi(x_{i})+(v_{\min,i}-x_{i})\cdot\nabla\psi(x_{i})+(\Delta x\odot S_{i}\odot z)\cdot\nabla\psi(x_{i})\\
	& = & \psi(x_{i})-\frac{\Delta x}{2}\cdot(S_{i}\odot\nabla\psi(x_{i}))+(\Delta x\odot S_{i}\odot z)\cdot\nabla\psi(x_{i})\\
	& = & \min\psi_{i}+(\Delta x\odot S_{i}\odot z)\cdot\nabla\psi(x_{i})\\
	& = & \min\psi_{i}+z\cdot\left(\Delta x\odot S_{i}\odot\nabla\psi(x_{i})\right)\\
	& = & \min\psi_{i}+z\cdot D_{i}
\end{eqnarray*}

\subsection{Partial volumes\label{sec:volumes}}

Our volumetric integration of $h\|\nabla\psi\|$ over the progressively
expanding $T$-sublevel set $\psi(x)\!<\!T$ will be carried out numerically
via cell-wise addition by approximating $ $$h\|\nabla\psi\|$ as
constant over each grid cell $x_{i}$ and then multiplying this by
the estimated volume of the cell subset where $\psi(x)\!<\!T$. We
will refer to volume within the cell of this $T$-sublevel set as
the \emph{partial volume} of the cell, denoted by $V_{i}(T)$ to reflect
its dependence on the level set value $T$. We will consider all of
the cell to belong to this $T$-sublevel set if $\max\psi_{i}\le T$
(thereby setting $V_{i}(T)$ to full cell volume $\prod\Delta x$);
we will consider a cell to be fully excluded if $\min\psi_{i}\ge T$
(thereby setting $V_{i}(T)$ to zero); and we will consider a cell
to be partially included and partially excluded whenever $\min\psi_{i}<T<\max\psi_{i}$.
Consistent with our convention to determine fully included and fully
excluded cells, we will use the same linear approximation $\psi_{i}$
within each cell to calculate $V_{i}(T)$ for partially included and
partially excluded cells as follows.

\begin{eqnarray*}
	V_{i}(T) & = & \int_{{\rm cell}(\Delta x,x_{i})}H\left(T-\psi_{i}(x)\right)dx\qquad\mbox{(where \ensuremath{H} denotes the Heaviside function)}\\
	& = & \int_{{\rm cell}(\Delta x,x_{i})}H\left(\frac{T-\psi_{i}(x)}{\|D_{i}\|_{1}}\right)dx\\
	& = & \left(\prod\Delta x\right)\underbrace{\int_{0}^{1}\cdots\int_{0}^{1}H\left(\frac{T-\min\psi_{i}-z\cdot D_{i}}{\|D_{i}\|_{1}}\right)dz_{1}\ldots dz_{n}}_{v}\\
	& = & \left(\prod\Delta x\right)\int_{0}^{1}\cdots\int_{0}^{1}H\left(\frac{T-\min\psi_{i}}{\max\psi_{i}-\min\psi_{i}}-z\cdot\frac{D_{i}}{\|D_{i}\|_{1}}\right)dz_{1}\ldots dz_{n}\\
	& = & \left(\prod\Delta x\right)v\left(\underbrace{\frac{T-\min\psi_{i}}{\max\psi_{i}-\min\psi_{i}}}_{\tau},\underbrace{\frac{D_{i}}{\|D_{i}\|_{1}}}_{d}\right)
\end{eqnarray*}
In this last line, we have defined the function $v$ which represents
the\emph{ fractional volume }of the unit cube given the \emph{fractional
	total increase} $\tau$ associated with the level set value $\tau$
and the vector $d$ of \emph{fractional directional increases} across
the cell. The partial volume $V_{i}$ is then obtained by multiplying
the total cell volume $\prod\Delta x$ by the fractional volume $v$
\[
v(\tau,d)=\int_{0}^{1}\cdots\int_{0}^{1}H\left(\tau-z\cdot d\right)\,dz_{1}\ldots dz_{n}
\]
Notice that if $\min\psi_{i}<T<\max\psi_{i}$ then 0<$\tau$<1. Moreover,
since the coefficients of $D_{i}$ are all non-negative, its $L^{1}$-norm
is given by the sum of its elements (the increases in $\psi_{i}$
across the cell along each orthogonal direction), meaning that the
elements $d_{1},\ldots,d_{n}$ of the $L^{1}$-normalized vector of
\emph{fractional directional increases} across the cell will sum to
1. In deriving dimension-specific formulas for $v$, it is helpful
to assume, without loss of generality (since the order of integration
with respect to the cell coordinates $z_{1},\ldots,z_{n}$ can be
permuted without consequence and since each cell coordinate ranges
across the same unit interval) that $d_{1},\ldots,d_{n}$ are sorted
in increasing order and therefore $d_{1}\le d_{2}\le\cdots\le d_{n}$.
When implementing the resulting formulas, however, it will be important
to sort the actual elements of $d$ accordingly before applying such
formulas.

\subsubsection{2D fractional volume (area) formulas}

\begin{figure}
	\centering{}%
	\begin{tabular}{ccc}
		\includegraphics{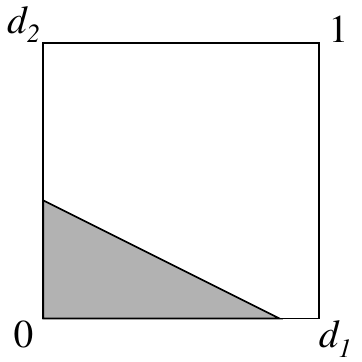} & \includegraphics{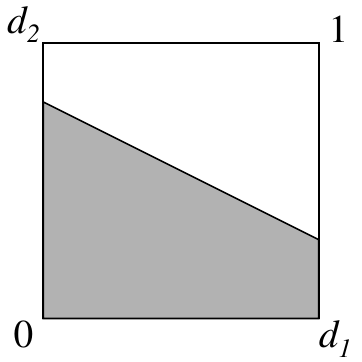} & \includegraphics{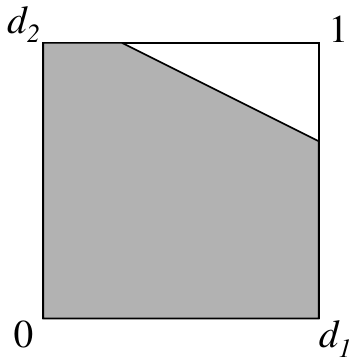}\tabularnewline
		$0<\tau<d_{1}$(case 1) & $d_{1}<\tau<d_{2}$ (case 2) & $d_{2}<\tau<1$ (case 3)\tabularnewline
	\end{tabular}\caption{\label{fig:2Dcases}Different possible fractional 2D cell geometries}
\end{figure}
Here we derive the fractional volume formula for $v$ in the 2D case
(and refer to it, accordingly, as the fractional area). There are
three different geometries to consider in the 2D case, which are illustrated
in~\ref{fig:2Dcases} where the lower-left corner represents
the origin of the unit square, with cell coordinates (0,0) and where
the upper right corner has cell coordinates (1,1). The horizontal
axis is associated with cell coordinate $z_{1}$, which means that
$v$(1,0)=$d_{1}$, while the vertical axis is associated with cell
coordinate $z_{2}$, which means that $v$(0,1)=$d_{2}$. In the case
that 0<$\tau$<$d_{1}$, the fractional area $v(\tau,d)$ is represented
in the left image of~\ref{fig:2Dcases} by the lower-left triangle
with base width $\frac{\tau}{d_{1}}$ and height $\frac{\tau}{d_{2}}$.
As $\tau$ increases and eventually falls within the range $d_{1}$<$\tau$<$d_{2}$,
the fractional area $v$ is represented in the middle image of~\ref{fig:2Dcases}
by a trapezoid with base width 1, left height $\frac{\tau}{d_{2}}$
and right height $\frac{\tau-d_{1}}{d_{2}}$. Finally as $\tau$ increases
further to fall within the final range $d_{2}$<$\tau$<1, the fractional
area $v$ is represented in the right image of~\ref{fig:2Dcases}
by the complement of the upper-right triangle with base width $\frac{1-\tau}{d_{1}}$
and height $\frac{1-\tau}{d_{2}}$. Thus, using the standard geometric
area formulas for the triangle and trapezoid and further noting that
the $L^{1}$-normalization of $d$ means that $d_{2}$=1-$d_{1}$,
it is easy to derive the following fractional area formulas in 2D.
\[
v_{2{\rm D}}(\tau,d)=\begin{cases}
\frac{1}{2}\frac{\tau^{2}}{d_{1}d_{2}}, & \mbox{[1] }\;\,0\le\tau<d_{1}\\
\\
\frac{1}{2}+\frac{\tau-\frac{1}{2}}{d_{2}}=\frac{\tau-\frac{d_{1}}{2}}{d_{2}}, & \mbox{[2] }\;\,d_{1}\le\tau\le1-d_{1}\\
\frac{\tau-\frac{1}{2}d_{1}}{d_{2}}\\
1-\frac{1}{2}\frac{(1-\tau)^{2}}{d_{1}d_{2}}, & \mbox{[3] }\;\,1-d_{1}<\tau\le1
\end{cases}
\]

In the special case where the $\nabla\psi(x_{i})$ is parallel to
one of the Cartesian grid axes, this will mean, by our ordering assumption
, that $d_{1}$=0 and $d_{2}$=1 and therefore only the middle case
will apply, with the simplified formula of $v_{2D}(\tau,d)=\tau$.
This simple expression may also be interpreted as the ``fractional
length'' $v_{1{\rm D}}(\tau)$ of a 1D cell which is simply identical
to the total fractional increase $\tau$.
\[
v_{2{\rm D}}\left(\tau,\begin{bmatrix}d_{1}\!=\!0\\
d_{2}
\end{bmatrix}\right)=v_{1{\rm D}}(\tau)=\tau
\]

\subsubsection{General fractional volume formulas}

While it is convenient in 2D, due to the simplicity and limited number
of geometries, to directly derive the fractional area formulas based
on the custom geometries associated to each of the three separate
cases, such an approach becomes considerably more complicated in three
dimensions, and completely impractical in even higher dimensions.
A more systematic method is to start by noting that for 0<$\tau$<$d_{1}$
the fractional partial volume corresponds to an $n$-dimensional simplex
whose vertices consist of an orthogonal corner at the cell origin
(0,$\ldots$,0) together with the intersections of a hyperplane with
each of the orthogonal cell coordinate axes. Since each of these axis
intersections occur within the unit interval of the associated cell
coordinate (by our assumption that 0<$\tau$<$d_{1}$ together with
our ordering assumption $d_{1}\le\cdots\le d_{n}$), the resulting
simplex is fully contained within the cell and therefore its volume
can be equated with the fractional volume $v(\tau,d)$. Noting that
the simplex intersects the $k$'th axis at coordinate value $z_{k}=\frac{\tau}{d_{k}}$
(by solving $\tau-z\cdot d=0$ when all the other cell coordinates
are set to zero), the simplex volume $v_{(0,\ldots,0)}(\tau,d)$ may
be expressed as follows.
\[
v_{(0,\ldots,0)}(\tau,d)=\frac{1}{n!}\prod_{k=1}^{n}\frac{\tau}{d_{k}}=\frac{\tau^{n}}{n!\prod d}
\]
As $\tau$ increases to enter the next range interval $d_{1}<\tau<d_{2}$
this simplex volume will exceed the fractional cell volume since part
of this simplex falls outside of the unit cell (the portion crossing
the $z_{1}$=1 hyperplane). However, the portion of this simplex which
falls outside of the cell is another simplex within the $z_{1}$>0
half-space with an orthogonal corner at (1,0,$\ldots$) whose edge
length along each orthogonal $z_{k}$ axis is given by $\frac{\tau-d_{1}}{d_{k}}$.
As such, its volume is easily expressed as well by
\[
v_{(1,0,\ldots,)}(\tau,d)=\frac{1}{n!}\prod_{k=1}^{n}\frac{\tau-d_{1}}{d_{k}}=\frac{(\tau-d_{1})^{n}}{n!\prod d}
\]
allowing us to express the fractional cell volume $v(\tau,d)$ as
the difference between the full simplex volume $v_{(0,\ldots,0)}(\tau,d)$
and this \emph{spillover simplex volume} $v_{(1,0,\ldots)}(\tau,d)$.
As $\tau$ further increases to fall within the interval $d_{2}<\tau<d_{3}$
then additional spillover occurs within the half-space $z_{2}$>1
hyperplane, which also takes the form of a simplex with an orthogonal
corner at (0,1,0,$\ldots$) with edge lengths given by $\frac{\tau-d_{2}}{d_{k}}$
along each $z_{k}$ axis and whose volume can therefore be expressed
as follows.
\[
v_{(0,1,0,\ldots,)}(\tau,d)=\frac{1}{n!}\prod_{k=1}^{n}\frac{\tau-d_{2}}{d_{k}}=\frac{(\tau-d_{2})^{n}}{n!\prod d}
\]
While the fractional volume $v$ can still be expressed as the full
simplex volume $v_{(0,\ldots,0)}$ minus the \emph{total spillover
	volume} special care must now be taken in computing the spillover
volume in terms of the individual \emph{spillover simplex volumes}
$v_{(1,0\ldots)}$ and $v_{(0,1,0,\ldots)}$ due to the fact that
the respective half-spaces given by $z_{1}$>1 and $z_{2}$>1 intersect
within the quarter-space where both $z_{1},z_{2}>1$. As such, the
half-space spillover simplices may also intersect (which will be the
case if $\tau>d_{1}+d_{2}$). If so, their intersection is yet another
simplex within the $z_{1,}z_{2}>1$ quarter-space with an orthogonal
corner at (1,1,0,$\ldots$) and edge lengths $\frac{\tau-d_{1}-d_{2}}{d_{k}}$
along each $z_{i}$ axis and with the following volume. 
\[
v_{(1,1,0,\ldots,)}(\tau,d)=\frac{1}{n!}\prod_{k=1}^{n}\frac{\tau-d_{1}-d_{2}}{d_{k}}=\frac{(\tau-d_{1}-d_{2})^{n}}{n!\prod d}
\]
Thus, summing the two half-space simplex volumes $v_{(1,0,\ldots)}$
and $v_{(0,1,0,\ldots)}$, $v_{(0,\ldots,0)}$ will yield the correct
total spillover volume so long as the two spillover simplices do not
intersect (which will be the case when $\tau<d_{1}+d_{2}$) but will
otherwise end up double-counting the quarter-space simplex volume
$v_{(1,1,0,\ldots)}$ which will therefore need to be subtracted in
order to compute the total spillover volume.

A systematic pattern emerges as $\tau$ continues to increase. One
starts by computing the total simplex volume $v_{(0,\ldots,0)}$.
Then, for each $k$ such that $\tau>d_{k}$, the corresponding half-space
simplex volume $v_{(\delta_{1k},\ldots,\delta_{nk})}$ (where $\delta$
denotes the Kronecker delta) is subtracted, then for each $j\ne k$
such that $\tau>d_{j}+d_{k}$, the corresponding quarter-space simplex
volume $v_{(\delta_{1j}\!+\!\delta_{1k},\ldots,\delta_{nj}\!+\!\delta_{nk})}$
is added to correct for having been over-subtracted by the two half-simplex
volumes indexed by $j$ and $k$. Then for each $j\ne k\ne m$ such
that $\tau>d_{j}+d_{k}+d_{m}$, the corresponding $\frac{1}{8}$-space
simplex volume $v_{(\delta_{1j}\!+\!\delta_{1k}\!+\!\delta_{1m},\ldots,\delta_{nj}\!+\!\delta_{nk}\!+\!\delta_{nm})}$
is subtracted to correct for a net over-addition by the three added
quarter-simplex volumes associated with each of the two-index combinations
$jk,$ $jm$, and $km$ taken from the triplet $jkm$ and the three
subtracted half-simplex volumes associated with each of the one-index
combinations $j$, $k$, and $m$ taken from the same triplet. Next,
for each $ $quadruplet of distinct indices such that $\tau$ exceeds
the sum of the respective elements of $d$, the corresponding $\frac{1}{16}$-space
simplex volume is added to correct for the net over-subtraction by
the four subtracted $\frac{1}{8}$-space simplex volumes associated
with each triple-index combination from the quadruplet, by the added
six quarter-simplex volumes associated with each double-index combination
from the quadruplet, and by the subtracted four half-simplex volumes
associated with each single-index combination from the quadruplet.
The process of alternating subtraction and addition of simplex volumes
continues until ($n$-1)-tuplets of indices have been accordingly
processed. Note that the generic formula for simplex volumes to be
added/subtracted at various stages in this process is as follows,
\[
v_{z_{{\rm vertex}}}=\frac{1}{n!}\prod_{k=1}^{n}\frac{\tau-d\cdot z_{{\rm vertex}}}{d_{k}}=\frac{(\tau-d\cdot z_{{\rm vertex}})^{n}}{n!\prod d}
\]
where $z_{{\rm vertex}}$ denotes the vertex cell coordinates (all
of which are 0 or 1) for the orthogonal corner of the corresponding
spillover simplex (for half-space simplices a single coordinate will
be 1, for $ $quarter-space simplices two coordinates will be 1, and
so on). 

The entire algorithm can be summarized succinctly as follows:
\begin{enumerate}
	\item Initialize the ``running fractional volume'' v to the volume of
	the full simplex with orthogonal corner at the cell origin as follows
	\[
	v=\frac{\tau^{n}}{n!\prod d}
	\]
	\item Now, looping from $m$=1 to $n$-1:
\end{enumerate}
\begin{itemize}
	\item For each $m$-tuplet of distinct (independent of order) indices $k_{1},\ldots,k_{m}$
	between 1 and $n$$ $,
	\item if $\tau>d_{k_{1}}+d_{k_{2}}+\cdots+d_{k_{m}}$ then update the running
	sum as follows
	\[
	v=v+(-1)^{m}\frac{\left(\tau-d_{k_{1}}-d_{k_{2}}-\cdots-d_{k_{m}}\right)^{n}}{n!\prod d}
	\]
\end{itemize}
This procedure must be modified if any of the elements of the fractional
directional increases vector $d$ are zero (since the common denominator
of all the simplex volumes would be zero in such cases). Fortunately,
the modification is extremely simple. Namely, a reduced dimension
vector $\hat{d}$ is created by removing the zero entries of $d$,
and the above procedure is then applied on $\hat{d}$ (also substituting
$n$ as well for the reduced dimension $\hat{n}$ of $\hat{d}$) for
the resulting set of lower-dimensional simplices.

\subsubsection{3D fractional volume formulas}

\begin{figure}
	\centering{}%
	\begin{tabular}{ccc}
		\includegraphics[width=0.3\textwidth]{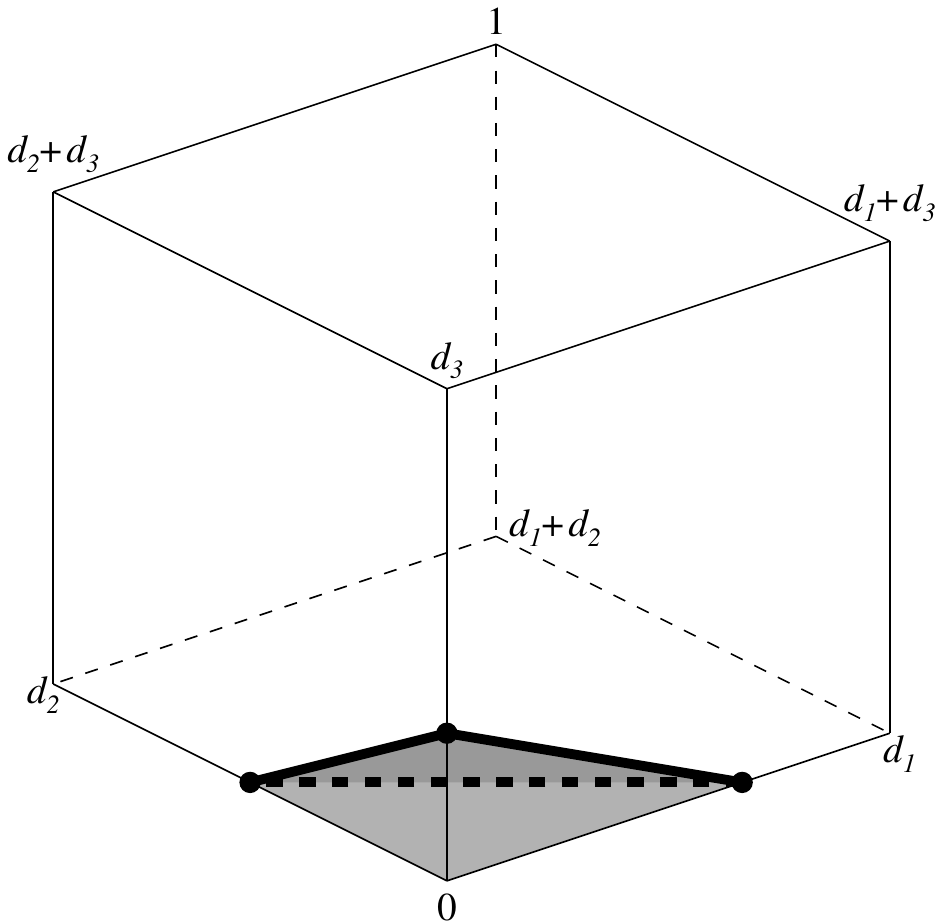} & \includegraphics[width=0.3\textwidth]{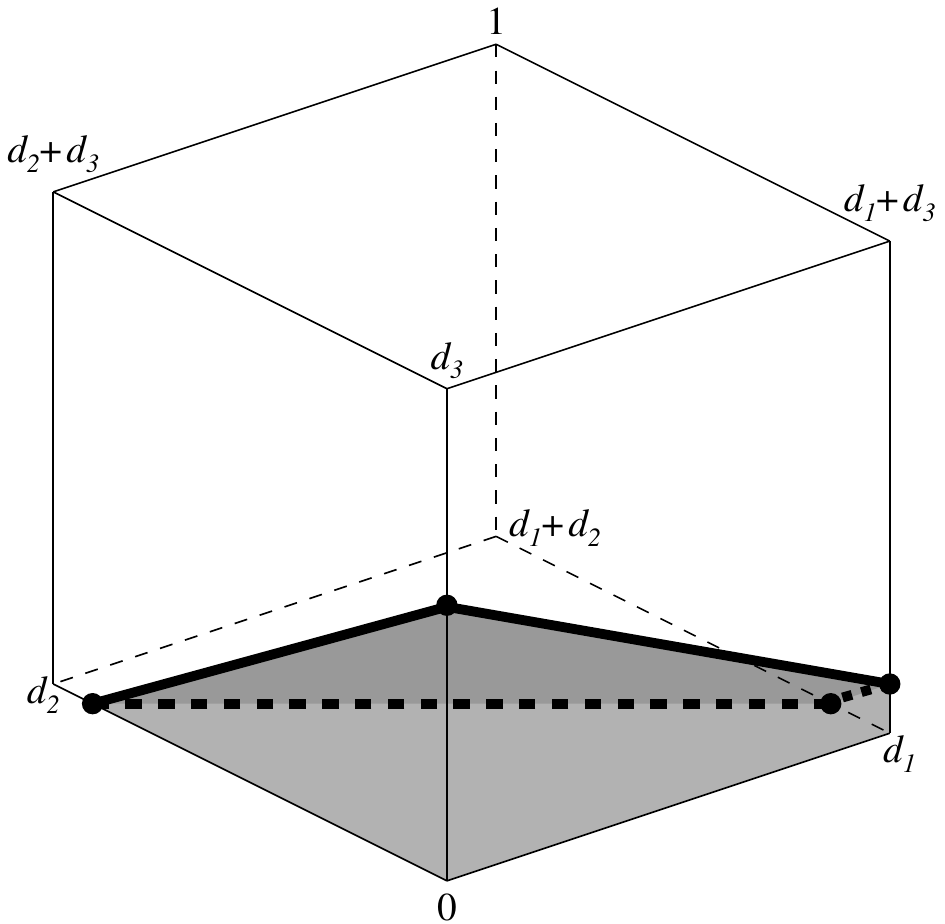} & \includegraphics[width=0.3\textwidth]{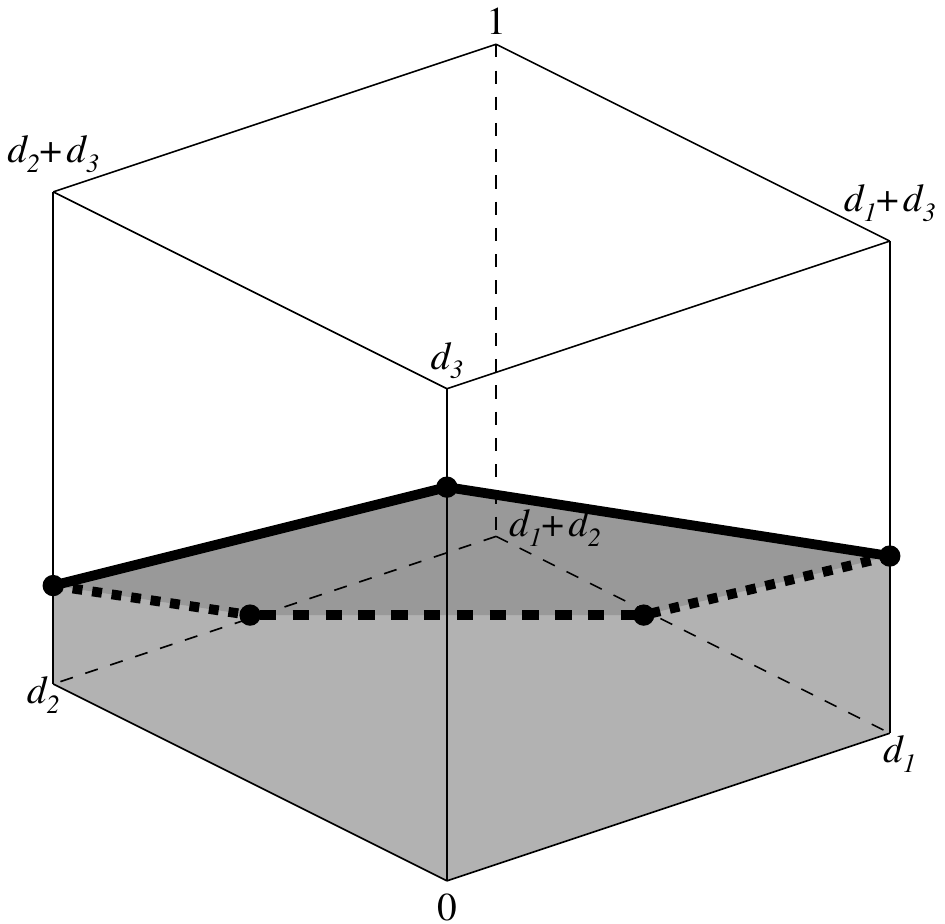}\tabularnewline
		$0<\tau<d_{1}$ (case 1) & $d_{1}<\tau<d_{2}$ (case 2) & $d_{2}<\tau<d_{3}$,\,$d_{1}+d_{2}$ (case 3)\tabularnewline
		\includegraphics[width=0.3\textwidth]{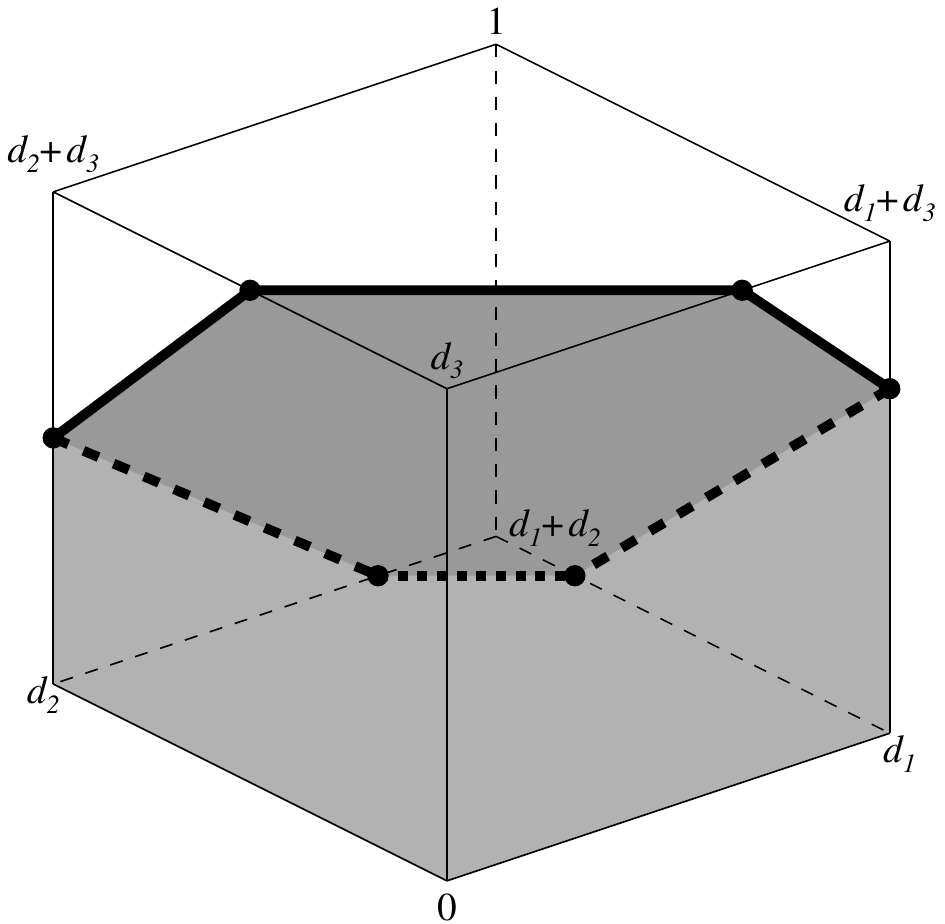} & \includegraphics[width=0.3\textwidth]{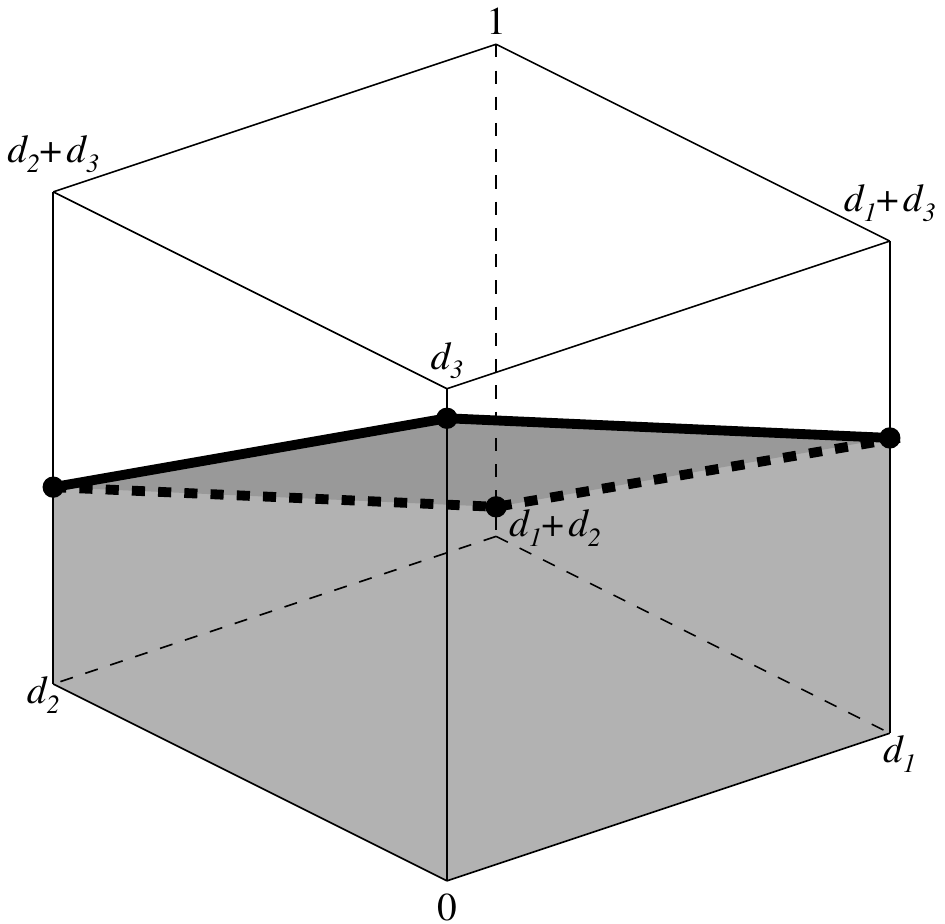} & \includegraphics[width=0.3\textwidth]{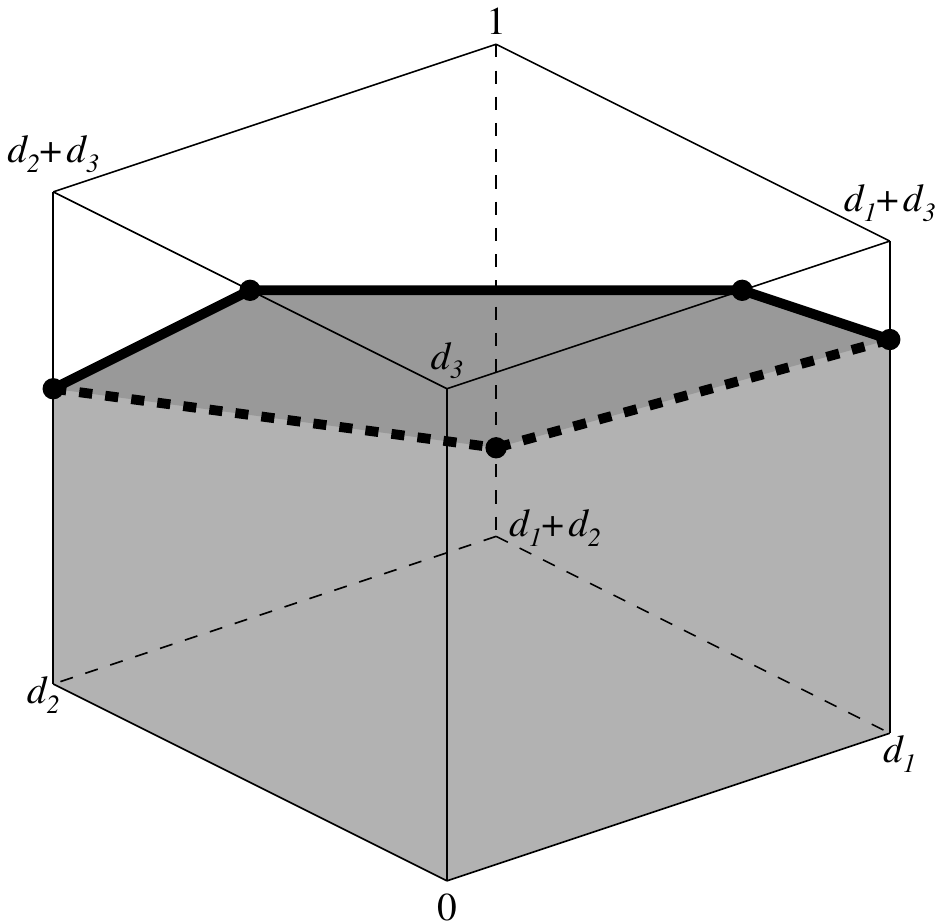}\tabularnewline
		$d_{3}<\tau<d_{1}+d_{2}$ (case 4a) & $d_{1}+d_{2}<\tau<d_{3}$ (case 4b) & $d_{1}+d_{2}$,\,$d_{3}<\tau<d_{1}+d_{3}$ (5) \tabularnewline
		\includegraphics[width=0.3\textwidth]{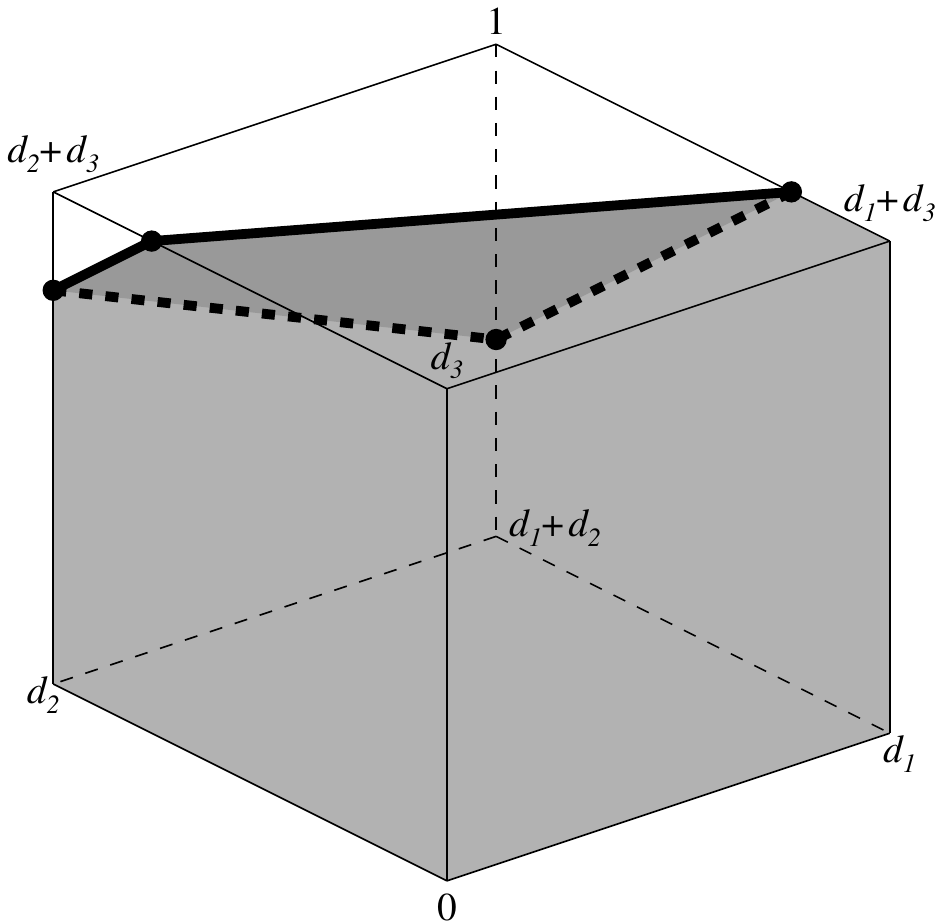} & \includegraphics[width=0.3\textwidth]{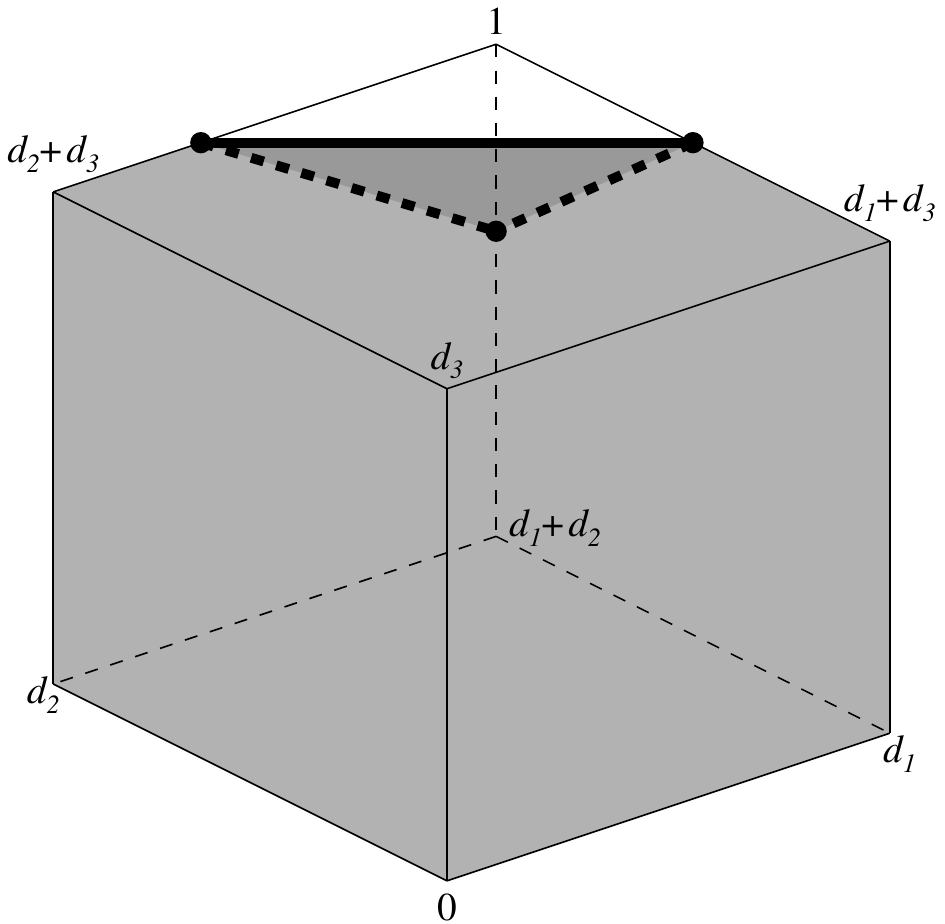} & \tabularnewline
		$d_{1}+d_{3}<\tau<d_{2}+d_{3}$ (case 6) & $d_{2}+d_{3}<\tau<1$ (case 7) & \tabularnewline
	\end{tabular}\caption{\label{fig:3Dcases}Different possible fractional 3D cell geometries}
\end{figure}

When applying this general strategy in 3D, it is helpful to note that
there are 6 different transition points, as the total fractional increase
$\tau$ progresses along the unit interval, across which the formulas
change. These transition values are given by $\tau=d_{1},d_{2},d_{3},d_{1}\!+\!d_{2},d_{1}\!+\!d_{3},d_{2}\!+\!d_{3}$,
which divide the interval into seven subintervals, but this only represents
one of two possible sorting orders. When sorting these transition
values, the first two values will always be $d_{1}$ and $d_{2}$
and the last two values will always be $d_{1}\!+\!d_{3}$ and $d_{2}\!+\!d_{3}$
based on our prior assumption that $d_{1}$, $d_{2}$, and $d_{3}$
are already sorted from smallest to largest to begin with. However,
this initial sorting assumption is not enough to tell us anything
about the relation between the middle two elements $d_{3}$ and $d_{1}\!+\!d_{2}$.
As such, we will consider separate scenarios where $d_{1}\!+\!d_{2}\le d_{3}$
and where $d_{1}\!+\!d_{2}>d_{3}$ (which is the same as $d_{3}\ge\frac{1}{2}$
and $d_{3}<\frac{1}{2}$ respectively since $d_{1}$+$d_{2}$+$d_{3}$=1),
thereby listing 8 possible subintervals even though only seven, at
most, are applicable for any given set of fractional increases $d_{1}$,
$d_{2}$, and $d_{3}$. 

Applying the strategy and notation developed for the general case
yields

\[
v_{3{\rm D}}=v_{(0,0,0)}-\begin{cases}
0, & \mbox{[1] }\;\,0<\tau<d_{1}\\
v_{(1,0,0)}, & \mbox{[2] }\;\,d_{1}<\tau<d_{2}\\
v_{(1,0,0)}+v_{(0,1,0)}, & \mbox{[3] }\;\,d_{2}<\tau<d_{3},\,d_{1}\!\!+\!d_{2}\\
v_{(1,0,0)}+v_{(0,1,0)}+v_{(0,0,1)}, & \mbox{[4a] }d_{3}<\tau<d_{1}\!\!+\!d_{2}\\
v_{(1,0,0)}+v_{(0,1,0)}\qquad\qquad-v_{(1,1,0)}, & \mbox{[4b] }d_{1}\!\!+\!d_{2}<\tau<d_{3}\\
v_{(1,0,0)}+v_{(0,1,0)}+v_{(0,0,1)}-v_{(1,1,0)}, & \mbox{[5] }\;\,d_{1}\!\!+\!d_{2},\,d_{3}<\tau<d_{1}\!\!+\!d_{3}\\
v_{(1,0,0)}+v_{(0,1,0)}+v_{(0,0,1)}-v_{(1,1,0)}\\-v_{(1,0,1)}, & \mbox{[6] }\;\,d_{1}\!\!+\!d_{3}<\tau<d_{2}\!+\!d_{3}\\
v_{(1,0,0)}+v_{(0,1,0)}+v_{(0,0,1)}-v_{(1,1,0)}\\-v_{(1,0,1)}-v_{(0,1,1)}, & \mbox{[7] }\;\,d_{2}\!+\!d_{3}<\tau<1
\end{cases}
\]
which, without any algebraic simplification (beyond factoring out
the common denominator for each simplex volume), results in the clearly
patterned formula
\[
\frac{1}{6d_{1}d_{2}d_{3}}\begin{cases}
\tau^{3}, & 0<\tau<d_{1}\\
\tau^{3}-(\tau\!-\!d_{1})^{3}, & d_{1}<\tau<d_{2}\\
\tau^{3}-(\tau\!-\!d_{1})^{3}-(\tau\!-\!d_{2})^{3}, & d_{2}<\tau<d_{3},\,d_{1}\!\!+\!d_{2}\\
\tau^{3}-(\tau\!-\!d_{1})^{3}-(\tau\!-\!d_{2})^{3}-(\tau\!-\!d_{3})^{3}, & d_{3}<\tau<d_{1}\!\!+\!d_{2}\\
\tau^{3}-(\tau\!-\!d_{1})^{3}-(\tau\!-\!d_{2})^{3}\\\qquad\qquad\quad+(\tau\!-\!d_{1}\!-\!d_{2})^{3}, & d_{1}\!\!+\!d_{2}<\tau<d_{3}\\
\tau^{3}-(\tau\!-\!d_{1})^{3}-(\tau\!-\!d_{2})^{3}-(\tau\!-\!d_{3})^{3}\\\qquad\qquad\quad+(\tau\!-\!d_{1}\!-\!d_{2})^{3}, & d_{1}\!\!+\!d_{2},\,d_{3}<\tau<d_{1}\!\!+\!d_{3}\\
\tau^{3}-(\tau\!-\!d_{1})^{3}-(\tau\!-\!d_{2})^{3}-(\tau\!-\!d_{3})^{3}\\+(\tau\!-\!d_{1}\!-\!d_{2})^{3}+(\tau\!-\!d_{1}\!-\!d_{3})^{3}, & d_{1}\!\!+\!d_{3}<\tau<d_{2}\!+\!d_{3}\\
\tau^{3}-(\tau\!-\!d_{1})^{3}-(\tau\!-\!d_{2})^{3}-(\tau\!-\!d_{3})^{3}\\+(\tau\!-\!d_{1}\!-\!d_{2})^{3}+(\tau\!-\!d_{1}\!-\!d_{3})^{3}+(\tau\!-\!d_{2}\!-\!d_{3})^{3}, & d_{2}\!+\!d_{3}<\tau<1
\end{cases}
\]
Despite the simplicity of this pattern (and the resulting convenience
in coding such case-by-case expressions through an accumulative loop
instead), better numerical precision can be obtained, in situations
where $d_{1}$ and possibly even $d_{2}$ are zero or close to zero,
by using the following case-by-case algebraically simplified expressions
(derived by exploiting $d_{1}$+$d_{2}$+$d_{3}$=1) despite the loss
of incremental looping convenience.
\[
v_{3{\rm D}}(\tau,d)=\begin{cases}
\frac{\tau^{3}}{6d_{1}d_{2}d_{3}}=\frac{1}{6}\,\frac{\tau}{d_{1}}\,\frac{\tau}{d_{2}}\,\frac{\tau}{d_{3}}, & \mbox{[1] }\;\,0\le\tau<d_{1}\\
\\
\frac{d_{1}^{2}+3\tau(\tau-d_{1})}{6d_{2}d_{3}}=\frac{1}{6}\,\frac{d_{1}}{d_{2}}\,\frac{d_{1}}{d_{3}}+\frac{1}{2}\,\frac{\tau-d_{1}}{d_{2}}\,\frac{\tau}{d_{3}}, & \mbox{[2] }\;\,d_{1}\le\tau<d_{2}\\
\\
\frac{1}{2}+\frac{\tau-\frac{1}{2}}{d_{3}}+\frac{(1-\tau-d_{3})^{3}}{6d_{1}d_{2}d_{3}}=\ldots\\\quad\frac{\tau-\frac{d_{1}+d_{2}}{2}}{d_{3}}+\frac{(1-\tau-d_{3})^{3}}{6d_{1}d_{2}d_{3}}, & \mbox{[3] }\;\,d_{2}\le\tau<d_{3},\,d_{1}\!+\!d_{2}\\
\\
\frac{1}{2}+\frac{\tau-\frac{1}{2}}{d_{3}}+\frac{(1-\tau-d_{3})^{3}-(\tau-d_{3})^{3}}{6d_{1}d_{2}d_{3}}=\ldots\\\quad\frac{\tau-\frac{d_{1}+d_{2}}{2}}{d_{3}}+\frac{(1-\tau-d_{3})^{3}-(\tau-d_{3})^{3}}{6d_{1}d_{2}d_{3}}, & \mbox{[4a] }d_{3}\le\tau\le d_{1}\!+\!d_{2}\quad\!{\scriptstyle \left[d_{3}\!<\!d_{1}\!+\!d_{2}\right]}\\
\\
\frac{1}{2}+\frac{\tau-\frac{1}{2}}{d_{3}}=\frac{\tau-\frac{d_{1}+d_{2}}{2}}{d_{3}}, & \mbox{[4b] }d_{1}\!+\!d_{2}\le\tau\le d_{3}\quad\!{\scriptstyle \left[d_{1}\!+\!d_{2}\!\le\!d_{3}\right]}\\
\\
\frac{1}{2}+\frac{\tau-\frac{1}{2}}{d_{3}}-\frac{(\tau-d_{3})^{3}}{6d_{1}d_{2}d_{3}}=\ldots\\\quad\frac{\tau-\frac{d_{1}+d_{2}}{2}}{d_{3}}-\frac{(\tau-d_{3})^{3}}{6d_{1}d_{2}d_{3}}, & \mbox{[5] }\;\,d_{1}\!+\!d_{2},\,d_{3}<\tau\le d_{1}\!+\!d_{3}\\
\\
1-\frac{d_{1}^{2}+3(1-\tau)(1-\tau-d_{1})}{6d_{2}d_{3}}=\ldots\\\quad1-\frac{1}{6}\,\frac{d_{1}}{d_{2}}\,\frac{d_{1}}{d_{3}}-\frac{1}{2}\,\frac{1-\tau-d_{1}}{d_{2}}\,\frac{1-\tau}{d_{3}}, & \mbox{[6] }\;\,d_{1}\!+\!d_{3}<\tau\le d_{2}\!+\!d_{3}\\
\\
1-\frac{(1-\tau)^{3}}{6d_{1}d_{2}d_{3}}=1-\frac{1}{6}\,\frac{1-\tau}{d_{1}}\,\frac{1-\tau}{d_{2}}\,\frac{1-\tau}{d_{3}}, & \mbox{[7] }\;\,d_{2}\!+\!d_{3}<\tau\le1
\end{cases}
\]
The above formulae are well conditioned numerical expressions which do not cause numerical underflow or overflow as the individual components tend to zero. If instead we were to stick to the earlier patterned formulae, they would behave badly as $d_1$ $d_2$ tend to zero. Notice that if $d_{1}$=0, then only cases 2, 4b, and 6 yield non-empty
subintervals, and the resulting three formulas match those of the
2D fractional area with $d_{2}$ and $d_{3}=1\!-\!d_{2}$ taking the
place of $d_{1}$ and $d_{2}=1\!-\!d_{1}$. If $d_{2}$=0 as well,
then $d_{3}$=1, and only the 4b subinterval is non-empty, and its
formula matches $\tau$ itself (the 1D fractional length).
\[
v_{3{\rm D}}\left(\tau,\begin{bmatrix}d_{1}\!=\!0\\
d_{2}\\
d_{3}
\end{bmatrix}\right)=v_{2{\rm D}}\left(\tau,\begin{bmatrix}d_{2}\\
d_{3}
\end{bmatrix}\right)\quad\mbox{and}
\]
\[\quad v_{3D}\left(\tau,\begin{bmatrix}d_{1}\!=\!0\\
d_{2}\!=\!0\\
d_{3}
\end{bmatrix}\right)=v_{2{\rm D}}\left(\tau,\begin{bmatrix}d_{2}\!=\!0\\
d_{3}
\end{bmatrix}\right)=v_{1{\rm D}}(\tau)
\]

\subsection[Efficient computation of M(T)]{Efficient computation of $M(T)$}

We now numerically approximate the volume integral $M(T)$ for each
of the discretized level set values $T=T_{\min},\,T_{\min}\!+\!\Delta T,\,T_{\min}\!+\!2\Delta T,\ldots,T_{\max}$
by summing over all cell-wise integrals $M_{i}(T)$ 

\begin{subequations}
	\begin{align}
	M(T) & \approx  \sum_{i}M_{i}(T)\label{eq:Msum}\\
	M_{i}(T) & = h(x_{i})\|\nabla\psi(x_{i})\|\underbrace{\left(\prod\Delta x\right)\,v\!\left(\frac{T-\min\psi_{i}}{\max\psi_{i}-\min\psi_{i}},\frac{D_{i}}{\|D_{i}\|_{1}}\right)}_{V_{i}(T)}\label{eq:Mi}
	\end{align}
\end{subequations}

computed according to the linear approximation $\psi_{i}$ of $\psi$
detailed in~\ref{sec:linear} and a constant approximation
$h_{i}$=$h(x_{i})$ of $h$ within each grid cell . While this mismatch
in approximation order (only constant for $h$ but linear for $\psi$)
may initially seem odd, it is justified by the fact that $M(T)$ is
defined by the volumetric integral of $h\|\nabla\psi\|$, and computation
of this integrand minimally requires a constant approximation of $h$
and a first-order approximation of $\psi$ (which results in a constant
approximation of the integrand $h\|\nabla\psi\|$ itself).

Rather than visiting all the grid cells every time we sum $M_{i}(T)$
for each discrete value of $T$, we may make the procedure more efficient
by partitioning the entire set of grid cells into subsets according
to whether their partial volumes $V_{i}(T)$ are ``full'' (equal
to the full cell volume), ``empty'' (equal to zero), or ``strictly
partial'' (non-zero but less then the full cell volume). This partition
will change as $T$ increases, but the number of ``full'' cells
can never decrease since such cells are characterized by $\max\psi_{i}\le T$,
while the number of ``empty'' cells can never increase since such
cells are characterized by $\min\psi_{i}\ge T$. Since we may represent
a single grid cell by its index $i$, we may accordingly represent
all of the cells grouped within these classes via the following index
sets
\begin{eqnarray*}
	{\cal J}_{{\rm full}}(T) & = & \{i\,|\,\max\psi_{i}\le T\}\\
	{\cal J}_{{\rm partial}}(T) & = & \{i\,|\,\min\psi_{i}<T<\max\psi_{i}\}\\
	{\cal J}_{{\rm empty}}(T) & = & \{i\,|\,\min\psi_{i}\ge T\}\\
	{\cal J}_{{\rm notfull}}(T) & = & \{i\,|\,\max\psi_{i}>T\}\quad=\quad{\cal J}_{{\rm partial}}(T)\cup{\cal J}_{{\rm empty}}(T)
\end{eqnarray*}
where we have added a fourth class ``not full'' which represents
cells that are either ``empty'' or ``partial''. 

The key point to exploit is that $M_{i}$ becomes constant for increasing
values of $T$ once a cell becomes full. Thus, rather than computing
$M_{i}(T+\Delta T)$ for any of the cells in ${\cal J}_{{\rm full}}(T)$,
the following relationship may be used instead
\begin{eqnarray*}
	\sum_{i}M_{i}(T+\Delta T) & = & \sum_{i\in{\cal J}_{{\rm full}}(T)}M_{i}(T+\Delta T)+\sum_{i\in{\cal J}_{{\rm notfull}}(t)}M_{i}(T+\Delta T)\\
	& = & \underbrace{\sum_{i\in{\cal J}_{{\rm full}}(T)}M_{i}(T)}_{M_{{\rm full}}(T)}+\sum_{i\in{\cal J}_{{\rm notfull}}(T)}M_{i}(T+\Delta T)
\end{eqnarray*}
where $M_{{\rm full}}(T)$ denotes the partial sum of $M_{i}(T)$
over only the cells that were already full at the level set value
$T$. So if we know $M_{{\rm full}}(T)$, then we do not need to visit
any of the ${\cal J}_{{\rm full}}(T)$ grid points to compute our
approximation of $M(T+\Delta T$) by using the formula
\begin{eqnarray*}
	M(T+\Delta T) & \approx & \sum_{i}M_{i}(T+\Delta T)\\
	& = & M_{{\rm full}}(T)+\sum_{i\in{\cal J}_{{\rm notfull}}(T)}M_{i}(T+\Delta T)
\end{eqnarray*}
which only requires we visit the ${\cal J}_{{\rm notfull}}(T)$ grid
points. Of course, in the process of visiting these ${\cal J}_{{\rm notfull}}(T)$
grid cells to compute their respective $M_{i}(T+\Delta T)$ contributions,
we may encounter grid cells that have now become full for $T+\Delta T$,
thereby reducing the size of the updated set ${\cal J}_{{\rm notfull}}(T+\Delta T)$
of grid cells that remain unfilled. The partial sum $M_{{\rm full}}$
should therefore be updated as well
\[
M_{{\rm full}}(T+\Delta T)=M_{{\rm full}}(T)+\!\!\!\!\underset{i\notin{\cal J}_{{\rm notfull}}(T\!+\!\Delta T)}{\sum_{i\in{\cal J}_{{\rm notfull}}(t)}}\!\!\!\!M_{i}(T+\Delta T)
\]
by adding the values of $M_{i}(T+\Delta T)$ for the newly filled
grid cells that were removed to create the updated $ $${\cal J}_{{\rm notfull}}(T+\Delta T)$
set.

\subsubsection{General traversal algorithm}\label{subsubsec:AlgoGeneral}

Applying these iterative concepts yields the following algorithm
\begin{description}
	\item [{Initialize}]~
\end{description}
\begin{enumerate}
	\item Set $M_{{\rm full}}$ initially to zero
	\item Put all grid indices in the initial list ${\cal J}_{{\rm notfull}}$
	\item For each grid index $i$, compute and store $\|\nabla\psi(x_{i})\|$,
	$\min\psi_{i}$, $\max\psi_{i}$, and the (sorted) relative increases
	$d$
\end{enumerate}
\begin{description}
	\item [{Loop}] $T=T_{\min}$ through $T=T_{\max}$ (by $\Delta T$ steps)
\end{description}
\begin{enumerate}
	\item Initialize $M(T)=M_{{\rm full}}$ 
	\item For each index $i$ within the list ${\cal J}_{{\rm notfull}}$
\end{enumerate}
\begin{itemize}
	\item compute $\Delta M=M_{i}(T)=\left(\prod\Delta x\right)h(x_{i})\|\nabla\psi(x_{i})\|\,v\!\left(\frac{T-\min\psi_{i}}{\max\psi_{i}-\min\psi_{i}},\frac{D_{i}}{\|D_{i}\|_{1}}\right)$
	\item increment $M(T)=M(T)+\Delta M$
	\item if $\max\psi_{i}\le T$ then also increment $M_{{\rm full}}=M_{{\rm full}}+\Delta M$
	and remove $i$ from the list ${\cal J}_{{\rm notfull}}$
\end{itemize}

\subsubsection{Even more efficient strategy for distance functions}\label{subsubsec:AlgoDistFunc}

In the case where $\psi$ represents a distance function to the $T_{\min}$
level set or, more generally, a function for which it is possible
to connect any domain point to the $T_{\min}$ level set by a trajectory
along which $\psi$ changes monotonically, an even more efficient
strategy can be developed. The idea here is to traverse the grid points
monotonically with respect to $\psi$, much like what is done in fast
marching algorithms (which are often employed to compute the level
set function $\psi$ itself). In this case, it is not necessary to
visit all of the previously ``empty cells'' when computing $M(T+\Delta T)$
since many of them will remain empty, and therefore contribute nothing
to the sum. While this is also typically true in the general case
outlined above, the lack of a causal monotonic structure makes it
impossible to know which of the previously empty cells for $t$ will
remain empty for $T+\Delta T$. However, when the causal structure
applies, we can be assured that the set of previously empty cells
which become nonempty for the new level set value of $T+\Delta T$
must be connected to the set of previously non-empty cells, thereby
allowing us to devise a more efficient search strategy as follows.
\begin{description}
	\item [{Initialize}]~
\end{description}
\begin{enumerate}
	\item Set $M_{{\rm full}}$ initially to zero
	\item For each grid index $i$, compute and store $\|\nabla\psi(x_{i})\|$,
	$\min\psi_{i}$, $\max\psi_{i}$, and the (sorted) relative increases
	$d$
	\item Put all grid indices for which $\min\psi_{i}<T_{\min}$ in the initial
	list ${\cal J}_{{\rm partial}}$
	\item Set a binary flag ${\rm visited_{i}}$ to one for each of these indices
	in the list, and set such a flag to zero for all remaining grid indices.
	We will refer to these two groups of indices as ``visited'' and
	``unvisited'' respectively.
\end{enumerate}
\begin{description}
	\item [{Loop}] $T=T_{\min}$ through $T=T_{\max}$ (by $\Delta T$ steps)
\end{description}
\begin{enumerate}
	\item Initialize $M(T)=M_{{\rm full}}$ 
	\item Retrieving each index $i$ sequentially from ${\cal J}_{{\rm partial}}$
	until the end of the list is reached
\end{enumerate}
\begin{itemize}
	\item compute $\Delta M=M_{i}(T)=\left(\prod\Delta x\right)h(x_{i})\|\nabla\psi(x_{i})\|\,v\!\left(\frac{T-\min\psi_{i}}{\max\psi_{i}-\min\psi_{i}},\frac{D_{i}}{\|D_{i}\|_{1}}\right)$
	\item increment $M(T)=M(T)+\Delta M$
	\item Append any ``unvisited neighbor index'' $j$ (meaning that $x_{j}$
	is a neighbor of $x_{i}$ in the Cartesian grid and that ${\rm visited_{j}}$=0)
	for which $\min\psi_{j}<T$ to the \emph{end} of the $ $list ${\cal J}_{{\rm partial}}$
	and then toggle the flag ${\rm visited_{j}}$=1.
	\item if $\max\psi_{i}\le T$ then also increment $M_{{\rm full}}=M_{{\rm full}}+\Delta M$
	and remove $i$ from the list ${\cal J}_{{\rm partial}}$ (and append
	any still unvisited neighbor index $j$ to ${\cal J}_{{\rm partial}}$
	and then toggle the flag ${\rm visited_{j}}$=1)
\end{itemize}
Note that the sequential traversal of the ${\cal J}_{{\rm partial}}$
list in step 2 is important since the list will potentially grow as
it is being traversed, with additional elements being appended to
the end of the list. If the list is not traversed sequentially (or
if appending does not occur at the end of the list) then the additional
elements many not be processed as they should be during the traversal.

\subsection[Direct computation of m(T)]{Direct computation of $m(T)$}

The level set integrals $m(T)$ of $h$ for each discretized level
set value $T$ are obtained by differentiating the antiderivative
function $M(T)$ which we have constructed thusfar. This may be done
numerically, by first computing the discrete representation of $M(T)$
at $T=T_{\min},\,T_{\min}\!\!+\!\Delta T,\ldots,\,T_{\max}$, and
then taking the numerical differences between consecutive samples
(and dividing by $\Delta T$). However, given that the approximation~\ref{eq:Msum} for $M$ consists of the sum of analytically differentiable
elements $M_{i}(T)$, one may instead approximate the function $m$
directly by differerentiating the elements $M_{i}$ inside the sum.
The sum of the resulting elements $m_{i}(T)$ may then be used to
directly approximate the function $m$ of level set integrals. 

\[
m(T)=\frac{d}{dT}M(T)\approx\frac{d}{dT}\sum_{i}M_{i}(T)=\sum_{i}m_{i}(T)
\]
These differentiated elements $m_{i}$ may be expressed as follows
\begin{eqnarray*}
	m_{i}(T) & = & \dot{M}_{i}(T)\\
	& = & \left(\prod\Delta x\right)h(x_{i})\|\nabla\psi(x_{i})\|\;\frac{d}{dT}\,v\!\left(\underbrace{\frac{T-\min\psi_{i}}{\max\psi_{i}-\min\psi_{i}}}_{\tau},\underbrace{\frac{D_{i}}{\|D_{i}\|_{1}}}_{d}\right)\\
	& = & \left(\prod\Delta x\right)h(x_{i})\frac{\|\nabla\psi(x_{i})\|}{\max\psi_{i}-\min\psi_{i}}\;\underbrace{\frac{\partial v}{\partial\tau}\!\left(\frac{T-\min\psi_{i}}{\max\psi_{i}-\min\psi_{i}},\frac{D_{i}}{\|D_{i}\|_{1}}\right)}_{\text{differential fractional volume}}
\end{eqnarray*}
where the partial derivative $\frac{\partial v}{\partial\tau}$ represents
the differential fractional volume within the unit cell. We obtain
this analytically by differentiating the expressions for the fractional
volume $v$ presented earlier (the 2D and 3D cases are presented below).
The surface area (arclength in 2D) element $dA(T)$ for the $T$-level
set of $\psi$ is therefore approximated within the cell by the intersected
area $A_{i}(T)$ of the $T$-level hypersurface of $\psi_{i}$ (intersected
length in 2D) given as follows,
\[
A_{i}(T)=\left(\prod\Delta x\right)\frac{\|\nabla\psi(x_{i})\|}{\max\psi_{i}-\min\psi_{i}}\;\frac{\partial v}{\partial\tau}\!\left(\frac{T-\min\psi_{i}}{\max\psi_{i}-\min\psi_{i}},\frac{D_{i}}{\|D_{i}\|_{1}}\right)
\]
yielding the following intuitive expression for the level set integral
cell-wise elements $m_{i}$.
\[
m_{i}(T)=h(x_{i})\,A_{i}(T)
\]
\subsubsection{2D differential fractional volume (area)}

The differential fractional volume $\frac{\partial v}{\partial\tau}$
in 2D (where it is actually better interpreted as a ``differential
fractional area'') is given as follows,

\[
\frac{\partial v_{2{\rm D}}}{\partial\tau}(\tau,d)=\frac{1}{d_{2}}\begin{cases}
\frac{\tau}{d_{1}}, & 0\le\tau<d_{1}\\
\\
1, & d_{1}\le\tau\le1-d_{1}\\
\\
\frac{1-\tau}{d_{1}}, & 1-d_{1}<\tau<1
\end{cases}
\]
with $\frac{\partial v_{2D}}{\partial\tau}$=0 for any other value
of $\tau$ outside the interval {[}0,1{]}.

\subsubsection{3D differential fractional volume}

The differential fractional volume $\frac{\partial v}{\partial\tau}$
in 3D is given as follows.

\[
\frac{\partial v_{3{\rm D}}}{\partial\tau}(\tau,d)=\frac{1}{d_{3}}\begin{cases}
\frac{\tau^{2}}{2d_{1}d_{2}}, & \mbox{[1] }\;\,0\le\tau<d_{1}\\
\\
\frac{\tau-\frac{1}{2}d_{1}}{d_{2}}, & \mbox{[2] }\;\,d_{1}\le\tau<d_{2}\\
\\
1-\frac{(1-\tau-d_{3})^{2}}{2d_{1}d_{2}}, & \mbox{[3] }\;\,d_{2}\le\tau<d_{3},\,d_{1}\!+\!d_{2}\\
\\
1-\frac{(1-\tau-d_{3})^{2}+(\tau-d_{3})^{2}}{2d_{1}d_{2}}, & \mbox{[4a] }d_{3}\le\tau\le d_{1}\!+\!d_{2}\quad{\scriptstyle \left[d_{3}\!<\!d_{1}\!+\!d_{2}\right]}\\
\\
1, & \mbox{[4b] }d_{1}\!+\!d_{2}\le\tau\le d_{3}\quad{\scriptstyle \left[d_{1}\!+\!d_{2}\!\le\!d_{3}\right]}\\
\\
1-\frac{(\tau-d_{3})^{2}}{2d_{1}d_{2}}, & \mbox{[5] }\;\,d_{1}\!+\!d_{2},\,d_{3}<\tau\le d_{1}\!+\!d_{3}\\
\\
\frac{(1-\tau)-\frac{1}{2}d_{1}}{d_{2}}, & \mbox{[6] }\;\,d_{1}\!_{-}\!d_{3}<\tau\le d_{2}\!+\!d_{3}\\
\\
\frac{(1-\tau)^{2}}{2d_{1}d_{2}}, & \mbox{[7] }\;\,d_{2}\!+\!d_{3}<\tau\le1
\end{cases}
\]
with $\frac{\partial v_{3D}}{\partial\tau}$=0 for any other value
of $\tau$ outside the interval {[}0,1{]}.

\section{Experimental results}
\label{sec:experiments}

In this section we present some numerical experimental results. Note that, our method is fundamentally designed to compute a whole family of integrals efficiently in a coupled way but nevertheless it is instructive to compare our accuracy to traditional single isosurface methods designed to perform integration along specific level surfaces. In a later sub-section we show more relevant experiments by computing integrals along families of isosurfaces. Here we compare our method to the Geometric numerical integration method of Min-Gibou~\cite{MinGibou07}. In~\cite{MinGibou07}, the authors had compared their Geometric numerical integration method to that of the first-order and second-order Delta function formulation~\cite{Smereka06}. Comparisons were made on accounts of robustness to perturbations of the interface location on the grid, order of convergence and numerical accuracy. We present our results together with that of~\cite{MinGibou07} and~\cite{Smereka06}, effectively comparing all four methods. Through experiments in both 2D and 3D domains, we demonstrate that our methods achieve similar levels of accuracy to traditional single isosurface methods while providing the advantage of being able to compute integrals along entire sets of level surfaces efficiently in a coupled way. We conduct all our experiments on an computer with a quad-core Intel core i7 CPU with 16 GB available RAM and we implemented our methods in the C++ programming language.

\subsection{Experiments in 2D}

In this section we compute the arc-length an ellipse represented as the zero level set of $\phi(x,y) = \frac{x^2}{1.5^2} + \frac{y^2}{0.75^2} - 1$ on a two dimensional grid spanning the range $[-5\: 5]$ in both dimensions. The true arc-length of this ellipse is $\approx 7.266336165$~\cite{Smereka06}. We used numerical implementation of our algorithm to experimentally compute the arc-length of the ellipse for various isotropic grid resolutions($\Delta x$) going from coarse resolution of $0.2$ to much finer grid resolution of $0.0625$. Similar to the experiments in~\cite{MinGibou07}, we conducted 50 trials where the ellipse was randomly translated on the grid. We generated a set of 50 pairs of uniform random numbers in the range of $[0,1)$ and used them as the center of the ellipse which effectively shifted the location of the interface on the grid. We used this same set of center perturbations for all grid resolutions($\Delta x)$.~\ref{tab:ellipseLength-50} shows the data for calculating the arc-length of the ellipse, averaged over 50 trials for each grid resolution. In order to compare our computation to the true arc-length of the ellipse we calculate the relative error for each trial and then compute the average relative error over the 50 trials for each grid resolution. We also compute additional statistics such as the minimum(Min) and maximum(Max) relative error and its standard deviation(SD). Finally, we also show the experimental order of convergence(Order) as we get finer grid resolutions and also show the ratio of the maximum and minimum relative error(Max/Min) over all 50 trials for a grid resolution. The order of convergence is calculated as $log2\left(\frac{E_k}{E_{k+1}}\right)/log_2\left(\frac{\Delta x ^k}{\Delta x^{k+1}}\right)$ where $E_k$ is the quantity under consideration for e.g. average relative error at a particular grid resolution level indicated by $k$.

\begin{table}[tbhp]
	{\footnotesize
		\caption{Computing arc-length of an ellipse for various grid resolutions for 50 trials}\label{tab:ellipseLength-50}
		\begin{center}
			\begin{tabular}{lcccccccc} \hline
				$\Delta x$ & Average & Order & SD & Min & Order & Max & Order & $\frac{max}{min}$\vspace{2 pt}\\  \hline
				\multicolumn{2}{l}{Intrinsic Integration}\\
				0.2&2.95E-03&0.00&2.11E-03&3.87E-06&0.00&1.04E-02&0.00&2684.60\\
				0.1&1.04E-03&1.51&5.84E-04&1.04E-04&-4.74&2.25E-03&2.21&21.71\\
				0.05&3.04E-04&1.77&1.80E-04&1.65E-05&2.65&6.24E-04&1.85&37.87\\
				0.025&8.90E-05&1.77&6.35E-05&3.00E-07&5.78&2.36E-04&1.40&786.26\\
				0.0125&3.12E-05&1.51&2.52E-05&8.54E-07&-1.51&1.20E-04&0.98&140.59\\
				0.00625&1.29E-06&1.28&9.40E-06&2.25E-07&1.93&2.99E-05&2.00&133.36\\
				& & & & & & & & \\
				\multicolumn{2}{l}{Min-Gibou~\cite{MinGibou07}}\\
				0.2&5.04E-03&0.00&2.15E-04&4.63E-03&0.00&5.49E-03&0.00&1.19\\
				0.1&1.26E-03&2.00&3.23E-05&1.17E-03&1.99&1.30E-03&2.08&1.11\\
				0.05&3.14E-04&2.00&6.61E-06&3.03E-04&1.95&3.26E-04&2.00&1.08\\
				0.025&7.84E-05&2.00&1.25E-06&7.50E-05&2.02&7.99E-05&2.03&1.07\\
				0.0125&1.96E-05&2.00&2.15E-07&1.90E-05&1.98&1.99E-05&2.01&1.04\\
				0.00625&4.90E-06&2.00&3.18E-08&4.83E-06&1.98&4.94E-06&2.01&1.02\\
				& & & & & & & & \\
				\multicolumn{4}{l}{First-order delta function approach~\cite{Smereka06}}\\
				0.2&8.96E-03&0.00&7.93E-03&1.28E-04&0.00&2.67E-02&0.00&208\\
				0.1&2.70E-03&1.73&2.96E-03&9.13E-05&1.31&1.07E-02&0.49&118\\
				0.05&9.55E-04&1.50&1.12E-03&4.19E-07&1.28&4.43E-03&7.77&10600\\
				0.025&3.21E-04&1.57&3.58E-04&7.32E-06&1.54&1.52E-03&-4.12&208\\
				0.0125&1.13E-04&1.51&1.22E-04&9.15E-06&1.53&5.28E-04&-0.32&57.7\\
				0.00625&3.94E-05&1.52&4.17E-05&2.01E-06&1.52&1.84E-04&2.19&91.7\\
				& & & & & & & & \\
				\multicolumn{4}{l}{Second-order delta function approach~\cite{Smereka06}}\\
				0.2&3.23E-03&0.00&2.74E-03&6.07E-04&0.00&1.30E-02&0.00&21.5\\
				0.1&5.74E-04&2.49&5.25E-04&2.93E-06&7.69&3.02E-03&2.10&1030\\
				0.05&1.13E-04&2.34&4.08E-05&2.55E-05&3.12&2.04E-04&3.88&8.01\\
				0.025&3.08E-05&1.87&8.31E-06&1.56E-05&0.70&4.72E-05&2.11&3.00\\
				0.0125&7.61E-06&2.01&1.51E-06&1.37E-06&3.50&1.21E-05&1.96&8.82\\
				0.00625&1.89E-06&2.01&1.82E-07&1.61E-06&0.23&2.12E-06&2.45&1.37		
			\end{tabular}
		\end{center}
	}
\end{table}
As we can see from~\ref{tab:ellipseLength-50}, the average relative error in calculating the arc-length of the ellipse is very similar to the Geometric Integration approach of Min-Gibou. On some accounts, such as minimum error, our approach performs better while in other cases such maximum-to-minimum ratio it is slightly worse. Overall, we see however the comparable accuracy to earlier single isosurface integration methods(although not the intended use of our methods).

\subsection{Experiments in 3D}
In this section we compute the surface area and volume of an ellipsoid represented as the zero level set of $\phi(x,y,z) = \frac{x^2}{1.5^2} + \frac{y^2}{0.75^2} + \frac{z^2}{0.5^2} - 1$. The true surface area of this ellipsoid is $\approx 9.901821$~\cite{Smereka06}. Similar to the experiments in~\cite{MinGibou07} and our experiments in 2D as shown in last sub-section, we conducted 50 trials where the ellipsoid was randomly translated on the grid. This time we generate a triplet of uniform random numbers in range $[0,1)$ to perturb the center of the ellisoid in order to shift the location of its interface on the grid. We again use this same set of 50 random translations for all grid resolutions and compute similar statistics related to the relative error as we did in the 2D case.~\ref{tab:ellipsoidArea-50} shows these statistics for calculating the surface area of the ellipsoid. As we can see our algorithm is much more robust to grid perturbations in 3D than in 2D as indicated by the low Max/Min errors (last column). It is more robust than the first and second order delta formulations and close to Min-Gibou method. The results are also close to second order accurate, closer to the method of Min-Gibou and better than the first and second order delta formulations. Finally, in~\ref{tab:ellipsoidVolume-50} we compute the volume of the same ellipsoid and compare the average relative errors to the Geometric Integration method. Again, we can see that our method is very competitive both in terms of accuracy and order of convergence. In this case, we conduct only one trial without any center perturbations.

\begin{table}[tbhp]
	{\footnotesize
		\caption{Computing surface area of an ellipsoid for various grid resolutions for 50 trials}\label{tab:ellipsoidArea-50}
		\begin{center}
			\begin{tabular}{lccccccccc} \hline
				$\Delta x$ & Average & Order & SD & Min & Order & Max & Order & $\frac{max}{min}$\vspace{2 pt}\\  \hline
				\multicolumn{2}{l}{Intrinsic Integration}\\
				0.2&2.15E-02&0.00&2.85E-03&1.57E-02&0.00&2.91E-02&0.00&1.85\\
				0.1&5.27E-03&2.03&1.04E-03&2.88E-03&2.45&7.75E-03&1.91&2.70\\
				0.05&1.36E-03&1.96&3.07E-04&4.93E-04&2.55&1.83E-03&2.09&3.71\\
				0.025&3.24E-04&2.07&4.23E-05&2.34E-04&1.07&4.38E-04&2.06&1.87\\
				& & & & & & & &\\
				\multicolumn{2}{l}{Min-Gibou~\cite{MinGibou07}}\\
				0.2&3.17E-02&0.00&2.90E-04&3.12E-02&0.00&3.22E-02&0.00&1.03\\
				0.1&7.91E-03&1.98&1.02E-05&7.89E-03&1.98&7.94E-03&2.02&1.00\\
				0.05&1.98E-03&2.00&6.81E-07&1.98E-03&2.00&1.98E-03&2.00&1.00\\
				0.025&4.94E-04&2.00&1.13E-07&4.94E-04&2.00&4.95E-04&2.00&1.00\\
				& & & & & & & &\\
				\multicolumn{4}{l}{First-order delta approach~\cite{Smereka06}}\\
				0.2&3.03E-02&0.00&7.12E-03&1.75E-02&0.00&4.73E-02&0.00&1.49\\
				0.1&7.77E-03&1.96&2.26E-03&3.95E-03&2.14&1.32E-03&1.84&1.51\\
				0.05&2.12E-03&1.87&7.36E-04&6.39E-04&2.62&4.48E-03&1.56&2.16\\
				0.025&5.20E-04&2.03&1.36E-04&3.41E-04&0.91&8.51E-04&2.39&1.54\\
				& & & & & & & &\\
				\multicolumn{4}{l}{Second-order delta approach~\cite{Smereka06}}\\
				0.2&6.86E-02&0.00&7.60E-02&1.01E-02&0.00&4.46E-01&0.00&44.4\\
				0.1&1.33E-02&2.37&8.38E-03&9.47E-04&3.41&2.91E-02&3.93&3.08\\
				0.05&2.60E-03&2.35&2.11E-03&2.49E-04&1.92&1.06E-02&1.45&42.4\\
				0.025&8.00E-04&1.70&1.01E-03&4.40E-06&5.82&6.88E-03&0.62&1570	
			\end{tabular}
		\end{center}
	}
\end{table}
\vspace{-6pt}
%

%
\begin{table}[tbhp]
	{\footnotesize
		\caption{Convergence rate for computing volume of ellipsoid for various grid resolutions for one trial. Relative error and order of convergence shown for our method and that of Min-Gibou}\label{tab:ellipsoidVolume-50}
		\begin{center}
			\begin{tabular}{lcccc} \hline
				$\Delta x$ & Intrinsinsic Integration & Order & Min-Gibou & order\\  \hline
				0.1&1.04E-02&0.00&1.36E-02&0.00\\
				0.05&2.59E-03&2.00&3.40E-03&2.00\\
				0.025&6.48E-04&2.00&8.50E-04&2.00				
			\end{tabular}
		\end{center}
	}
\end{table}

%
%
\subsection{Computing families of integrals}\label{subsec:WholeFamilyIntegrals}
In the last two sub-sections we presented some results and comparisons for computing integrals along a single levelset of a function. We did that with the aim of demonstrating that our method achieves accuracy similar to traditional single isocontour methods even though developed with aim of computing a whole family of integrals in one shot in contrast to these traditional methods. In this sub-section we show that our techniques allows simultaneous and coupled computation of families of integrals while offering significant computational advantages.

In this experiment, we compute the arc length of whole families of level curves of a circle represented as the zero level set of $\phi(x,y) = x^2 + y^2 - 1$ on a grid spanning the range $[-5\; 5]$ in both X and Y directions with an isotropic grid resolution $\Delta x = 0.025$. We compute the arc lengths of progressively expanding circles corresponding to level curve represented by $T_{min} = -0.5$ to $T_{max} = 2.0$. We divide the range $[T_{min}, T_{max}]$ into $nTs$ equally spaced intervals with $nTs$ ranging from 5 to 200 in increments of 5. We then employ the two algorithms detailed in~\ref{subsubsec:AlgoGeneral} and~\ref{subsubsec:AlgoDistFunc}. We refer to the method of~\ref{subsubsec:AlgoGeneral} as ``Coupled Non-Causal Integration'' and the method of~\ref{subsubsec:AlgoDistFunc} which applies to the special case of distance functions as ``Coupled Causal Integration''. We compare the accuracy and speed of our methods with a brute-force approach of using the standard marching triangles based explicit isosurfacing technique in a sequential manner, one at a time over all the desired sub-level surfaces.

\begin{figure}[ht]
	\centering
	\includegraphics[width=\textwidth]{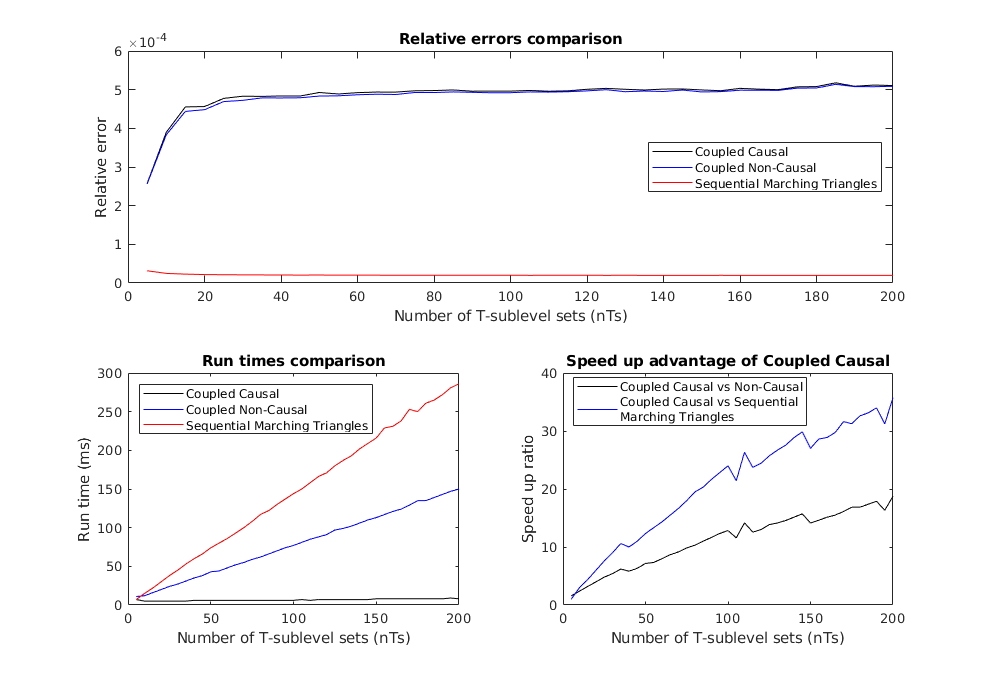}
	\vspace{-15pt}
	\caption{Accuracy(top) and run times(bottom) comparison of integration methods based on the experimental setup described in~\ref{subsec:WholeFamilyIntegrals}. The accuracy (relative error) of sequential marching triangle is approximately $10^{-05}$ on average compared to our two methods at $10^{-04}$. In terms of run times comparison(bottom left), our methods specially designed to compute integrals along families of T sub-level curves/surfaces significantly outperform the application of standard marching triangles sequentially to each of the $nTs$ sub-level curves. In terms of computational speed up(bottom right), in the special case of integration along signed distance functions and when the number of T sub-level curves are very large, our Coupled Causal method(\ref{subsubsec:AlgoDistFunc}) delivers a speed up of approximately 15 times compared to our Coupled Non-Causal method(\ref{subsubsec:AlgoGeneral}) and approximately 35 times compared to sequential application of standard marching triangle approach to all the $nTs$ sub-level curves.}
	\label{fig:intrinsicExperiments}
\end{figure}

In~\ref{fig:intrinsicExperiments} top row,  we compare the relative errors in computing the arc-length of expanding concentric circles using our two methods and marching triangles based explicit isosurfacing scheme. We take the average of the relative errors for each sub level curve $\tau \in [T_{min}, T_{max}]$. We can see that the accuracy of marching triangles is higher at the order of approximately $10^{-05}$. However, our two methods are also not far behind at approximately $10^{-04}$. The main advantage is in the computational efficiency of our methods compared to traditional methods. In~\ref{fig:intrinsicExperiments} bottom left, we can see that the run time of both our methods is consistently lower than marching triangles based approach. The run times of the Coupled Causal method which applies to distance functions is significantly lower than either of the other methods. Finally, in~\ref{fig:intrinsicExperiments} bottom right, we can see that as the number of T-level surfaces over which integration is performed increases, the Coupled Causal method can provide significant computational speed up while maintaining relatively good accuracy. Compared to the Coupled Non-Causal traversal algorithm, it can speed up the calculation of integrals by 15 - 20 times and almost 25 - 35 times when compared to Marching Triangles based approach. The Coupled Causal algorithm runs in almost constant time even as the number of integrals grows.

\subsection{Practical Application}\label{subsec:PracticalApplication}
Identifying the precise location and extent of an object in an image is a fundamental task in many computer vision applications. Prior knowledge about the objects of interest may improve segmentation performance in the presence of noise, occlusion, and model errors. A review of these \textit{recognition}-segmentation approaches with different priors, such as color, motion, shape, and texture is given in \cite{Cremers2007}. Shape is a powerful feature for many applications, particularly for medical imaging where shapes of body parts do not vary much among patients. Shape models are plentiful but most relevant to our work are level set representations of shape, which have been incorporated as shape priors for active contour techniques, e.g. in~\cite{Bresson,  leventon, tsaiyezzi}.

Additional robustness can be achieved when combining shape with appearance. Examples of level set approaches are \cite{Fritscher20073DIS, Huangdeformable, Yangdeformable} where principal component analysis (PCA) is performed on the level sets and on the pixel-intensity image of the training set to obtain coupled shape and appearance models. These approaches typically require computationally expensive two-dimensional warps of intensity templates to migrate between different shape configurations. Utilizing the computationally efficient intrinsic integration approach developed in this paper, we can propose a smaller, transformed set of intensity features (a one-dimensional function) by numerically integrating image intensities along iso-contours of the object’s shape.

\begin{figure}[ht]
    \centering
    \begin{subfigure}{0.45\textwidth}
        \centering
        \includegraphics[width=6cm]{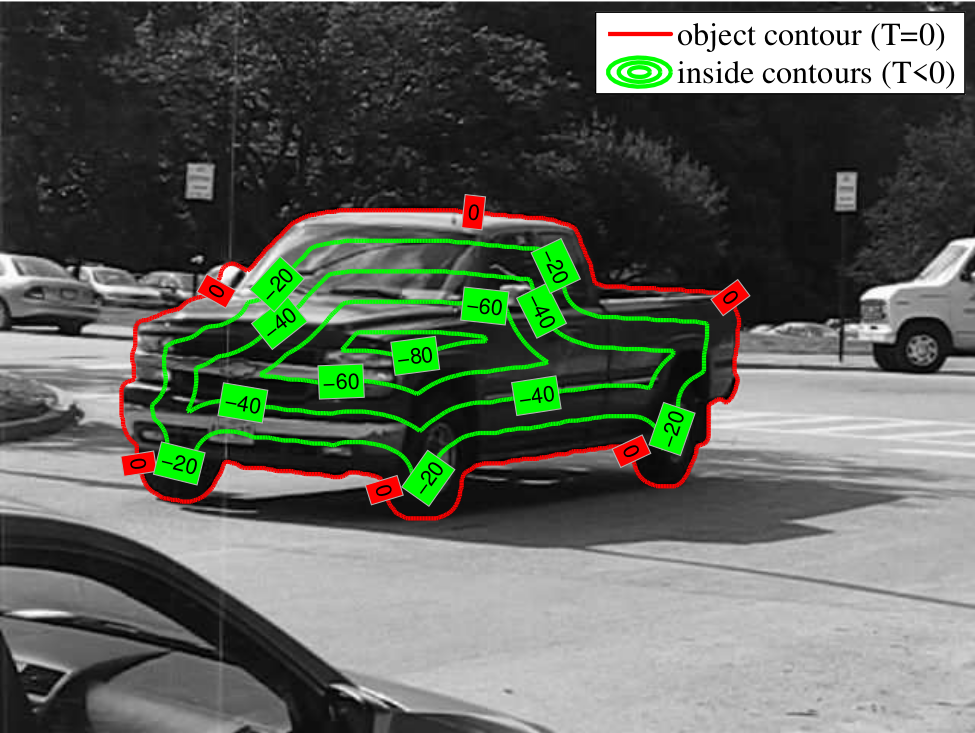}
        \caption{level set curves $C_T$}
	\label{fig:f_T_levelsets}
    \end{subfigure}%
    ~ 
    \begin{subfigure}{0.45\textwidth}
        \centering
        \includegraphics[height=2.5cm]{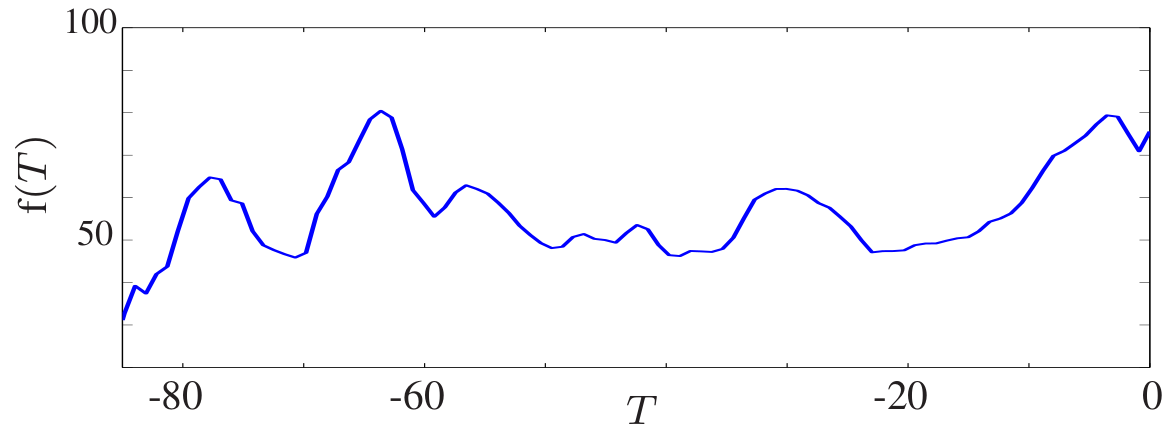}
        \caption{mean intensities along $T$-level-sets}
	\label{fig:fT_fT}
    \end{subfigure}
    \caption{Illustration of the photo-geometric representation.}
\label{fig:f_T}
\end{figure}

Given an image $I$ and the object boundary shape $C$ (and its equivalent signed-distance function $\psi$), traditional appearance models build a two-dimensional image template. Instead, we propose to use the mean image intensities along iso-contours of the object shape, i.e., our one-dimensional template is defined as 
\begin{equation}
\text{f}(T) = \frac{\int_{C_T}I\, ds} {\text{length}(C_T)},
\label{eq:fT2}
\end{equation}
where $\tf: [\psi_\text{min},\, 0] \mapsto \R$ are the mean intensities and $C_T = \{x: \psi(x)=T \}$ are the iso-contours. In \autoref{fig:f_T}, these quantities are illustrated for the image of a truck from the Berkeley Motion Segmentation dataset~\cite{brox}. The level sets $C_T$ for some values of $T$ between $\psi_\text{min}$ and $0$ are shown in \autoref{fig:f_T_levelsets} and the mean intensity along these curves is plotted in \autoref{fig:fT_fT}. 

We call $\text{f}(T)$ the coupled \textit{photo-geometric} representation of an object because it couples the object's geometric information $C_T$ with its photometric information $I$. As a feature, it is well-suited for training purposes because it is invariant to translation and rotation. It is also invariant to scale if the domain of $\tf(T)$ is normalized to a constant interval, e.g., $[-1,\, 0]$. As an appearance model, the photo-geometric representation can be viewed as a compromise between the powerful but cumbersome Active appearance models~\cite{cootes2001} and efficient but less discriminative methods, e.g., Chan-Vese~\cite{chanveseedges} where object intensity is modeled through a finite set of statistics. Our proposed algorithm using the one-dimensional photo-geometric descriptor is more general than finite sets of statistics, and at the same time not as cumbersome as two-dimensional template models.

A key point to note is that in order to compute this coupled \textit{photo-geometric} feature representation, we need to integrate over a number of level sets of an input image shape. There can be a multitude of training images and computing these features using traditional techniques applied sequentially to each level set will be extremely computationally expensive. In this case, we can apply our Coupled-Causal intrinsic integration technique (sec. \ref{subsec:WholeFamilyIntegrals}) to compute the whole family of required integrals for each image in one shot. See~\cite{mueller2021efficiently} for the complete development and implementation of a PCA based coupled shape and appearance model using the coupled \textit{photo-geometric} feature descriptor, efficiently computed using our integration techniques.

\section{Conclusions}
\label{sec:conclusions}

In this paper, we have presented a method of computing an entire continuous family of integrals over a whole family of T-level surfaces of a function by making use of the coarea formula. We have presented several experiments to show that our method is still comparable to traditional single isosurface integration methods in terms of accuracy achieved while making possible the task of computing whole family of integrals in a simultaneous efficiently coupled way. We additionally presented a special method applicable to integration over distance functions which provides significant speed up over other traditional single isocontour/surface methods while still maintaining comparable accuracy. Our methods can find applications over many problem domains but particularly Active Contour(Surface) based methods routinely used in Computer Vision problems.
\clearpage

\bibliography{references}

\end{document}